\documentclass[12pt]{article}
\usepackage[utf8]{inputenc}
\usepackage{amssymb,amsmath,amsfonts,eucal,mathrsfs,amsthm} 
\usepackage[colorinlistoftodos]{todonotes}
\usepackage[allcolors=blue]{hyperref}
\usepackage{empheq}
\newtheorem{theorem}{Theorem}
\newtheorem{proposition}[theorem]{Proposition}
\newtheorem{lemma}[theorem]{Lemma}

\newtheorem{corollary}[theorem]{Corollary}

\newtheorem{remarks}[theorem]{Remarks}

\theoremstyle{definition}
\newcommand{\R}{\mathbb{R}}
\newcommand{\Q}{\mathbb{Q}}
\newcommand{\Sf}{\mathbb{S}}

\newcommand{\spa}{\mbox{span}}
\newcommand{\hess}{\mbox{Hess\,}}
\newcommand{\Ric}{\mbox{Ric }}

\newcommand{\rank}{\mbox{rank }}

\newcommand{\nap}{\nabla^{\perp}}
\newcommand{\nab}{\tilde\nabla}

\newcommand{\End}{\mbox{End}}
\newcommand{\Hom}{\mbox{Hom}}
\newcommand{\trace}{\mbox{tr\,}}

\newcommand{\Y}{\mathcal{Y}\,}

\newcommand{\Les}{\mathbb{L}}
\def\<{{\langle}}
\def\>{{\rangle}}

\def\F{{\cal F}}

\def\T{{\cal T}}

\def\Y{{\cal Y}}
\def\n{\nabla}

\def\a{\alpha}

\def\be{\begin{equation} }
\def\ee{\end{equation} }

\def\proof{\noindent{\it Proof:  }}
\def\qed{\ifhmode\unskip\nobreak\fi\ifmmode\ifinner
\else\hskip5 pt \fi\fi\hbox{\hskip5 pt \vrule width4 pt
height6 pt  depth1.5 pt \hskip 1pt }}
\makeatletter
\newcommand{\subjclass}[2][]{\let\@oldtitle\@title
\gdef\@title{\@oldtitle\footnotetext{#1 
\emph{Mathematics Subject Classification:} #2}}}
\newcommand{\keywords}[1]{\let\@@oldtitle\@title
\gdef\@title{\@@oldtitle\footnotetext
{\emph{Key words and phrases.} #1.}}}
\makeatother
\begin{document}

\title{Infinitesimally  Moebius bendable hypersurfaces}
\maketitle
\begin{center}
\author{M. I. Jimenez        and
        R. Tojeiro$^*$}  
        \footnote{Corresponding author}
      \footnote{This research was initiated while the first author was supported by CAPES-PNPD Grant 88887.469213/2019-00 and was finished under the support of Fapesp Grant 2022/05321-9. The second author was partially supported by Fapesp grant 2022/16097-2 and CNPq grant 307016/2021-8.\\
      Data availability statement: Not applicable.}
\end{center}
\date{}

\begin{abstract}
   Li, Ma and Wang have provided in \cite{LMW} a partial classification of the so-called Moebius deformable hypersurfaces, that is, the umbilic-free Euclidean hypersurfaces $f\colon M^n\to \mathbb{R}^{n+1}$  that admit non-trivial deformations preserving the Moebius metric. For $n\geq 5$, the classification was completed by the authors in \cite{JT2}. In this article we obtain an infinitesimal version of that classification. Namely, we introduce the notion of an infinitesimal Moebius variation of an umbilic-free immersion $f\colon M^n\to \mathbb{R}^m$ into Euclidean space as a one-parameter family of immersions $f_t\colon M^n\to \mathbb{R}^m$, with $t\in (-\epsilon, \epsilon)$ and $f_0=f$, such that the Moebius metrics determined by $f_t$ coincide up to the first order. Then we characterize isometric immersions $f\colon M^n\to \R^m$  of arbitrary codimension that admit a non-trivial infinitesimal Moebius variation among those that admit a non-trivial conformal infinitesimal variation, and use such characterization to classify the umbilic-free Euclidean hypersurfaces of dimension $n\geq 5$ that admit non-trivial infinitesimal Moebius variations.
\end{abstract}

\vspace{2ex}


\date{}
\maketitle

\section{Introduction}
Given an  isometric immersion $f\colon M^n\to\R^{m}$ of a Riemannian manifold $(M^n,g)$ into Euclidean space with normal bundle-valued second fundamental form $\a\in \Gamma(\Hom(TM,TM;N_fM))$, let $\phi\in C^\infty(M)$ be defined by
\be\label{phi}
\phi^2=\frac{n}{n-1}(\|\a\|^2-n\|\mathcal{H}\|^2),
\ee
where $\mathcal{H}$ is the mean curvature vector field of $f$ and $\|\a\|^2\in C^\infty(M)$ is given at any point $x\in M^n$ by
$$
\|\a(x)\|^2=\sum_{i,j=1}^n\|\a(x)(X_i,X_j)\|^2,
$$
in terms of  an orthonormal basis $\{X_i\}_{1\leq i\leq n}$  of $T_xM$. Notice that $\phi$ vanishes precisely at the umbilical points of $f$.
The metric
$$
g^*=\phi^2 g,
$$
defined on the open subset of non-umbilical points of $f$, is a Moebius invariant metric called the \emph{Moebius metric} determined by $f$. Namely, if $\tilde f=\tau\circ f$ for some  Moebius transformation of $\R^m$, then the Moebius metrics of $f$ and $\tilde f$  coincide.

 It is a  fundamental fact, proved by Wang in \cite{Wa},  that a hypersurface  $f\colon M^n\to\R^{n+1}$ is uniquely determined, up to Moebius transformations of the ambient space, by its Moebius metric and its \emph{Moebius shape operator} $S=\phi^{-1}(A-HI)$, where $A$ is the shape operator of $f$ with respect to a unit normal vector field $N$ and $H$ is the corresponding mean curvature function. A similar result holds for submanifolds of arbitrary codimension  (see \cite{Wa} and Section $9.8$ of \cite{DT}).

Li, Ma and Wang have provided in \cite{LMW} a partial classification of the  hypersurfaces $f\colon M^n\to\R^{n+1}$, $n\geq 4$, that are not determined, up to Moebius transformations of $\R^{n+1}$,  only by their Moebius metrics, called \emph{Moebius deformable hypersurfaces}. 
For $n\geq 5$, the classification of Moebius deformable hypersurfaces was completed by the authors in \cite{JT2}. 

Moebius deformable hypersurfaces  belong to the more general class of \emph{conformally deformable} hypersurfaces, that is, the hypersurfaces  $f\colon M^n\to\R^{n+1}$ for which  $M^n$ admits an immersion  $\tilde f\colon M^n\to\R^{n+1}$  such that $f$ and $\tilde f$ induce conformal metrics on $M^n$ and do not differ by a Moebius transformation.  The study of conformally deformable hypersurfaces goes back to Cartan \cite{Ca2} (see also \cite{DT1} and Chapter $17$ of \cite{DT}).

  Our main goal in this article is to classify the \emph{infinitesimally} Moebius bendable hypersurfaces, that is, the umbilic-free hypersurfaces $f\colon M^n\to \mathbb{R}^{n+1}$ into Euclidean space that admit a one-parameter family of immersions $f_t\colon M^n\to \mathbb{R}^{n+1}$, with $t\in (-\epsilon, \epsilon)$ and $f_0=f$, whose Moebius metrics coincide with that of $f$ \emph{up to the first order}, in a sense that is made precise below.

   Let $f\colon M^n\to\R^m$ be an isometric immersion free of umbilical points.
We call a  smooth map $F\colon (-\epsilon, \epsilon)\times M^n\to\R^m$  a \emph{Moebius variation} of $f$ if   $f_t=F(t,\cdot)$, with $f_0=f$, is an immersion that determines the same Moebius metric for any $t\in (-\epsilon, \epsilon)$. In other words, if $g_t$ is the metric induced by $f_t$, then
$
g_t^*=g_0^*
$
for all $t\in I$. 

Trivial Moebius  variations can be produced by composing $f$ with the elements of a smooth one-parameter family of Moebius transformations of the Euclidean ambient space. Thus, the results in  \cite{JT2} and \cite{LMW} give a classification of the umbilic-free hypersurfaces  $f\colon M^n\to\R^{n+1}$  of dimension $n\geq 5$ that admit non-trivial Moebius  variations.

  We are interested in the  umbilic-free isometric immersions  $f\colon M^n\to\R^m$  that satisfy the weaker condition of admitting  non-trivial \emph{infinitesimal} Moebius  variations.  By an \emph{infinitesimal Moebius  variation} of an isometric immersion $f\colon M^n\to\R^m$ without umbilical points we mean a smooth map $F\colon (-\epsilon, \epsilon)\times M^n\to\R^m$ such that the maps $f_t=F(t,\cdot)$, with $f_0=f$, are immersions whose corresponding Moebius metrics coincide \emph{up to the first order}. This means that $\frac{\partial}{\partial t}|_{t=0}g_t^*=0$, that is,
$$
\frac{\partial}{\partial t}|_{t=0}(\phi_t^2\<f_{t*}X,f_{t*}Y\>)=0
$$
for all $X,Y\in\mathfrak{X}(M)$, where $\phi_t^2$ is given by \eqref{phi} for the immersion $f_t$, $t\in (-\epsilon, \epsilon)$. 

   Given a smooth variation  $F\colon (-\epsilon, \epsilon)\times M^n\to\R^m$ of  an isometric immersion $f\colon M^n\to\R^m$, one defines its \emph{variational vector field} by $\T = F_*\partial/\partial t|_{t=0}$. When the immersions  $f_t=F(t,\cdot)$ are the compositions of $f$ with the elements of a smooth one-parameter family of Moebius transformations of $\R^m$, the  variational vector field $\T$ is the restriction to $M^n$ of a conformal Killing vector field
of $\R^m$.  Accordingly, an infinitesimal Moebius  variation $F\colon (-\epsilon, \epsilon)\times M^n\to\R^m$ of an isometric immersion $f\colon M^n\to\R^m$ without umbilical points is said to be \emph{trivial} if the variational vector field $\T$ associated with $F$ is the restriction to $M^n$ of a conformal Killing vector field of $\R^m$. We say that $f$ is \emph{infinitesimally Moebius bendable} if it admits an infinitesimal  Moebius variation that is non-trivial restricted to any open subset of $M^n$. It is \emph{locally infinitesimally Moebius bendable} if each point $x\in M^n$ has an open neighborhood $U$ such that $f|_U$ is infinitesimally Moebius bendable.

  In order to state our classification of the umbilic-free infinitesimally Moebius bendable Euclidean hypersurfaces of dimension $n\geq 5$, we need some further definitions.  
 
 First, by a \emph{conformally surface-like hypersurface} $f\colon M^n\to \R^{n+1}$ we mean a hypersurface that differs by a Moebius transformation of $\mathbb{R}^{n+1}$ from either a cylinder or a rotation hypersurface over a surface in $\R^3$, or from a cylinder over a three-dimensional 
hypersurface of $\R^4$ that is a cone over a surface in $\Sf^3$. We say, accordingly, that $f$ is a conformally surface-like hypersurface \emph{determined by a surface $h\colon L^2\to \mathbb{Q}_{\epsilon}^3$}, with $\epsilon=0, -1$ or $1$, respectively. 

 Now we recall how a two-parameter family of hyperspheres  in $\R^{n+1}$ is determined by a surface $s\colon L^2\to \mathbb{S}_{1,1}^{n+2}$ into the Lorentzian sphere 
$$
\Sf_{1,1}^{n+2}=\{x\in\Les^{n+3}\colon\<x,x\>=1\}
$$
 in the Lorentz space $\mathbb{L}^{n+3}$.  
 
 Let $f\colon M^n\to \R^{n+1}$ be an oriented hypersurface with respect to a unit normal vector field~$N$. 
Then the family of hyperspheres 
$
x\in M^n\mapsto S(h(x),r(x))
$
with radius $r(x)$ and center 
$
h(x)=f(x)+r(x)N(x)
$
is enveloped by $f$.  If, in particular, $1/r$ is the mean curvature of $f$, it is called the \emph{central sphere congruence} of $f$.

Let $\mathbb{V}^{n+2}$ denote the light cone in $\mathbb{L}^{n+3}$ and let $\Psi=\Psi_{v,w,C}\colon\mathbb{R}^{n+1}\to\mathbb{L}^{n+3}$
be the isometric embedding onto 
$$
\mathbb{E}^{n+1}=\mathbb{E}^{n+1}_w=\{u\in\mathbb{V}^{n+2}:\<u,w\> =1\}\subset\mathbb{L}^{n+3}
$$
given by 
$$
\Psi(x)=v+Cx-\frac{1}{2}\|x\|^2w,
$$
in terms of $w\in \mathbb{V}^{n+2}$, $v\in \mathbb{E}^{n+1}$ 
and a linear isometry $C\colon\mathbb{R}^{n+1}\to\{v,w\}^\perp$. Then the congruence of hyperspheres  $x\in M^n\mapsto S(h(x),r(x))$ is determined by  
the map 
$S\colon M^n\to\mathbb{S}_{1,1}^{n+2}$ defined by
$$
S(x)=\frac{1}{r(x)}\Psi(h(x))+\frac{r(x)}{2}w,
$$
for $\Psi(S(h(x),r(x)))=\mathbb{E}^{n+1}\cap S(x)^\perp$ for all $x\in M^n$. 
The map $S$ has rank $0<k<n$, that is, it corresponds to a $k$-parameter congruence of hyperespheres, if and only if $\lambda=1/r$ is a principal curvature of $f$ with  constant multiplicity $n-k$ (see Section $9.3$ of \cite{DT} for details). In this case, $S$ gives rise to a map $s\colon L^k\to \mathbb{S}_{1,1}^{n+2}$ such that $S\circ \pi=s$, where $\pi\colon M^n\to L^k$ is the canonical projection onto the quotient space of leaves of $\ker (A-\lambda I)$.

    Finally, a surface $h\colon L^2\to \mathbb{Q}^3_\epsilon$ is said to be a \emph{generalized cone}
over  a unit-speed curve $\gamma\colon I\to \mathbb{Q}^2_c$, $c\geq \epsilon$, in an umbilical surface $\mathbb{Q}^2_c\subset \mathbb{Q}^3_\epsilon$, if
$L^2=I\times J$ is a product of intervals $I, J\subset \R$ and 
$h(s,t)=\exp_{\gamma(s)}tN(s)$
for all $(s,t)\in I\times J$, where $\exp$ is the exponential map of $\mathbb{Q}^3_\epsilon$ and $N$ is a unit normal vector field to $\mathbb{Q}^2_c$ along $\gamma$. Notice that $h$ has $0$ as one of its principal curvatures, with the $t$-coordinate curves as the correspondent curvature lines. Generalized cones without totally geodesic points are precisely the isothermic surfaces that have $0$ as a simple principal curvature. Recall that a surface $h\colon L^2\to \mathbb{Q}^3_\epsilon$ is
\emph{isothermic} if each non-umbilic point of $L^2$ has an open neighborhood where one can define isothermic (that is, conformal) coordinates whose coordinate curves are lines of curvature of $h$.

\begin{theorem} \label{thm:main}
Let $f\colon M^n\to \mathbb{R}^{n+1}$, $n\geq 5$, be an umbilic-free  infinitesimally Moebius bendable hypersurface.  Then there exists an open and dense subset $\mathcal{U}^*$ of $M^n$ such that $f$ is of one of the following types on each connected component $U$ of  $\mathcal{U}^*$:
\begin{itemize} 
\item[(i)] a conformally surface-like hypersurface determined by an isothermic surface $h\colon L^2\to \Q_\epsilon^3$, $\epsilon \in \{-1, 0, 1\}$.  
\item[(ii)] a hypersurface whose central sphere congruence is
determined by a  minimal space-like surface $s\colon L^2\to \Sf_{1,1}^{n+2}$.
\end{itemize}

In particular, $f$ has  a principal curvature with multiplicity $n-1$ or $n-2$ at any point of $M^n$, and the first possibility occurs on a connected component $U$ of  $\mathcal{U}^*$ if and only if $f$ is given on $U$ as in part $(i)$, with the surface  $h\colon L^2\to \Q_\epsilon^3$ being a generalized cone over a unit-speed curve $\gamma\colon J\to \Q_{c}^2$ in an umbilical surface $\Q_{c}^2\subset \Q_{\epsilon}^3$, $c\geq \epsilon$.

Conversely, any  simply connected hypersurface as in $(ii)$ is  infinitesimally Moebius bendable,  and for any hypersurface as in $(i)$ there exists an open dense subset where $f$ is locally infinitesimally Moebius bendable. 
\end{theorem}

It follows from Theorem \ref{thm:main} and the main result in \cite{JT2} that, within the class of hypersurfaces that  are not conformally surface-like on  any open subset and have a principal curvature of constant multiplicity $n-2$, the families of those that are either Moebius deformable or infinitesimally Moebius bendable coincide. On the other hand, among conformally surface-like hypersurfaces, the class of infinitesimally Moebius bendable hypersurfaces is strictly larger than that of Moebius deformable hypersurfaces. Indeed, while a surface in the former class  is determined by an arbitrary isothermic surface, the elements
in the latter are determined by particular isothermic surfaces, namely, Bonnet surfaces admitting isometric deformations preserving the mean curvature. 

   Our approach to prove Theorem \eqref{thm:main} is rather different from those used in both \cite{JT2}  or \cite{LMW} to classify the Moebius deformable hypersurfaces. It is based on the theory developed in \cite{DJ} and \cite{DJV} of the more general notions of \emph{conformal  variations} and \emph{conformal infinitesimal variations} of an isometric immersion $f\colon M^n\to\R^m$, which are natural generalizations of the corresponding classical concepts of isometric variations and isometric infinitesimal variations.
   
   A smooth map $F\colon (-\epsilon, \epsilon)\times M^n\to\R^m$
   is a \emph{conformal variation} of an isometric immersion  $f\colon M^n\to \R^m$ if the maps $f_t=F(t,\cdot)$, with $f_0=f$, are conformal immersions for any $t\in (-\epsilon, \epsilon)$, that is, if there is a positive  $\gamma\in  C^\infty((-\epsilon, \epsilon)\times M^n)$, with $\gamma(0, x) = 1$ for all $x\in M^n$, such that
\begin{equation}
\label{eq:confvar}\gamma(t, x)\<{f_t}_*X, {f_t}_*X\>=\<X,Y\>
\end{equation}
for all  $X, Y\in \mathfrak{X}(M)$, where $\<\, ,\,\>$ stands for the metrics of both $\mathbb{R}^m$ and $M^n$. Thus Moebius variations are particular cases of conformal variations for which $\gamma(t,x)=\phi_0^{-2}(x)\phi^2_t(x)$ for all $(t,x)\in I\times M^n$.

   \emph{Conformal infinitesimal variations} of an isometric immersion  $f\colon M^n\to \R^m$  are smooth variations for which \eqref{eq:confvar} holds up to the first order, that is,
\begin{equation}\label{inf}
\frac{\partial}{\partial t}|_{t=0} (\gamma(t, x)\<{f_t}_*X, {f_t}_*X\>)=0
\end{equation}
 for all  $X, Y\in \mathfrak{X}(M)$. Eq. \eqref{inf} implies that 
the variational vector field $\T = F_*\partial/\partial t|_{t=0}$ of $F$ satisfies
\begin{equation}\label{bending}
\<\tilde\nabla_X\T, f_*Y\>+\<f_*X, \tilde\nabla_Y\T\>=2\rho\<X,Y\>
\end{equation}
for all  $X, Y\in \mathfrak{X}(M)$, where $\rho(x)=-(1/2)\partial \gamma/\partial t(0,x)$.
For this reason, a smooth section $\T\in \Gamma(f^*T\mathbb{R}^m)$ that satisfies
\eqref{bending} is called a \emph{conformal infinitesimal bending} of $f$ with conformal factor $\rho\in C^{\infty}(M)$. In particular,  the variational vector field $\T=F_*\partial/\partial t|_{t=0}$ of an infinitesimal Moebius variation, which we call an \emph{an infinitesimal  Moebius bending},  is also a conformal infinitesimal bending of $f$ whose conformal factor is 
\begin{equation}\label{rhomoeb}
\rho=-\frac{1}{2}\frac{\partial}{\partial t}|_{t=0}(\gamma(t,x))=-\frac{1}{2}\phi_0^{-2}(x)\frac{\partial}{\partial t}|_{t=0}(\phi_t^2(x)).
\end{equation}

By the above, the variational vector  field of a conformal infinitesimal variation is a conformal infinitesimal bending. Conversely, any
conformal infinitesimal bending $\T\in \Gamma(f^*T\mathbb{R}^m)$ is the variational vector  field of a (non-unique) conformal infinitesimal variation $F\colon (-\epsilon, \epsilon)\times M^n\to \mathbb{R}^m$ of $f$. For instance, one may take
$$
{F}(t, x) = f(x) + t\T(x)
$$
for all $(t,x)\in (-\epsilon, \epsilon)\times M^n$.  The reason why it is convenient to consider the conformal infinitesimal bending associated with a conformal infinitesimal variation of an isometric immersion  $f\colon M^n\to \R^m$  is that one can establish a fundamental theorem providing necessary and sufficient conditions for the existence of a  conformal infinitesimal bending (and hence of a conformal infinitesimal variation); see \cite{DJ}.

  \emph{Infinitesimal variations} of an isometric immersion $f\colon M^n\to \R^m$ correspond to the conformal infinitesimal variations for which the function $\gamma$ in \eqref{inf} has the constant value $\gamma=1$. The associated variational vector fields are called  \emph{infinitesimal bendings} and correspond to the conformal infinitesimal bendings with conformal factor $\rho=0$. The reason why isothermic surfaces $f\colon L^2\to \Q_c^3$ appear in this context is that they are precisely the surfaces that are locally \emph{infinitesimally Bonnet bendable}, that is, the surfaces that admit local infinitesimal variations $F\colon (-\epsilon, \epsilon)\times L^2\to\Q_c^3$ such that the mean curvature functions $H_t$ of $f_t=F(t,\cdot)$, $t\in (-\epsilon, \epsilon)$, coincide up  to the first order, that is,  
$\partial/\partial t|_{t=0}{H}_t=0$ (see, e.g., Proposition $9$ of \cite{JT}).
  
  The study of  hypersurfaces $f\colon M^n\to\R^{n+1}$, $n\geq 3$, that admit non-trivial variations preserving the \emph{induced} metric goes back to Sbrana \cite{sb} and Cartan~\cite{ca1} (see also \cite{DFT2} or Chapter $11$ of \cite{DT}), whereas 
the  hypersurfaces that admit non-trivial infinitesimal variations were investigated by Sbrana \cite{sb0} (see also \cite{DV}, Chapter $14$ of \cite{DT} and \cite{DJ1}). We point out that the latter class turns out to be much larger than the former.
   
In the proof of Theorem \eqref{thm:main},  a main step is the following characterization of independent interest of the  infinitesimally Moebius bendable isometric immersions $f\colon M^n\to \R^m$  of arbitrary codimension among those that admit a non-trivial conformal infinitesimal bending $\T$ with conformal factor $\rho\in C^{\infty}(M)$. In the next statement, we denote by $\mathcal{H}$ the mean curvature vector field of $f$ and by $\mathcal{L}\in\Gamma(N_fM)$ the normal vector field given by
\be\label{meanbeta}
\mathcal{L}=\frac{1}{n}\sum_{i=1}^n\beta(X_i,X_i)\in \Gamma(N_fM)
\ee
for any orthonormal frame $\{X_1, \ldots, X_n\}$ of $M^n$, where $\beta$ is the symmetric section of $\mbox{Hom}\,(TM,TM; N_fM)$ associated with  $\T$ (see \eqref{eq:beta} below).

  \begin{theorem}\label{prop:cibxmib}
An  isometric immersion $f \colon M^n\to \R^m$ is infinitesimally Moebius bendable if and only if it admits a non-trivial  conformal infinitesimal  bending 
such that
\be\label{eqrhomoeb}
\Delta\rho+n\<\mathcal{L},\mathcal{H}\>=0.
\ee
\end{theorem}

 By means of Theorem \ref{prop:cibxmib}, it is shown in the proof of Theorem \ref{thm:main} that any  conformal infinitesimal  bending of a hypersurface as in part $(ii)$ of the statement of that result is also an infinitesimal  Moebius bending.

\section{The Fundamental theorem of conformal infinitesimal bendings}

 In this section we  recall from \cite{DJ} the Fundamental theorem for conformal infinitesimal bendings of Euclidean hypersurfaces.

Let $f\colon M^n\to\R^m$ be an isometric immersion and let $\T$ be a conformal infinitesimal bending of $f$ with conformal factor $\rho$, that is, $\T$ and $\rho$ satisfy \eqref{bending}. 
Defining $L\in\Gamma(\Hom(TM;f^*T\R^m))$ 
by
$$
LX=\nab_X \T-\rho f_*X=\T_*X-\rho f_*X
$$
for any $X\in\mathfrak{X}(M)$, then \eqref{bending} can be written as 
$$
\<LX,f_*Y\>+\<f_*X,LY\>=0
$$
for all $X,Y\in\mathfrak{X}(M)$. Let 
$B\in \Gamma(\mbox{Hom}\,(TM,TM; f^*T\R^m))$ be  given by
$$
B(X,Y)=(\nab_XL)Y=\nab_XLY-L\n_XY
$$
for all $X,Y\in\mathfrak{X}(M)$, and define 
$\beta\in \Gamma(\mbox{Hom}\,(TM,TM; N_fM))$ by
\begin{equation}\label{eq:beta}
\beta(X,Y)=(B(X,Y))_{N_fM}=(\nab_X\nab_Y\T-\nab_{\nabla_XY}\T)_{N_fM}-\rho\a(X,Y)
\end{equation}
for all $X,Y\in\mathfrak{X}(M)$.  Flatness of the ambient 
space and the symmetry of $\alpha$ imply that $\beta$ is symmetric. 

Given  $\eta\in\Gamma(N_fM)$, let
 $B_\eta\in\Gamma(\End(TM))$ be given by 
$\<B_\eta X,Y\>=\<\beta(X,Y),\eta\>$ for all $X, Y\in \mathfrak{X}(M)$. Then it can be shown that
\begin{equation}\label{derGaussC}
A_{\beta(Y,Z)}X+B_{\alpha(Y,Z)}X-A_{\beta(X,Z)}Y-B_{\alpha(X,Z)}Y
+\<Y,Z\>\nabla_X\nabla\rho
\end{equation}
for all $X,Y,Z\in\mathfrak{X}(M)$;  see \cite{DJ}, where a fundamental theorem for  conformal infinitesimal bendings of Euclidean submanifolds with arbitrary codimension was obtained. 
Here we restrict ourselves to state that theorem for the particular case of hypersurfaces. 

Given a  hypersurface $f\colon M^n\to \mathbb{R}^{n+1}$, let  $A$ be its shape operator with respect to a unit normal vector field $N$, and let  $\mathcal{B}\in \Gamma(\End(TM))$ be given by 
$$\<\mathcal{B}X,Y\>=\<\beta(X,Y),N\>$$
for all $X, Y\in \mathfrak{X}(M)$. The Fundamental theorem for conformal infinitesimal bendings of $f$ reads as follows. 

\begin{theorem}\label{fundhyp} \emph{(\cite{DJ})}
The pair $(\mathcal{B}, \rho)$ associated with a conformal 
infinitesimal bending of the hypersurface  $f\colon M^n\to \mathbb{R}^{n+1}$ satisfies the equations
\be\label{hyp1}
\mathcal{B}X\wedge AY-\mathcal{B}Y\wedge AX+X\wedge \mbox{Hess}\,\rho(Y)-Y\wedge \mbox{Hess}\,\rho(X)=0
\ee
and
\be\label{hyp2}
(\nabla_X\mathcal{B})Y-(\nabla_Y\mathcal{B})X+(X\wedge Y)A\nabla\rho=0
\ee
for all $X,Y\in\mathfrak{X}(M)$. Conversely, if  $M^n$ is simply connected, then a symmetric tensor  $\mathcal{B}\in \Gamma(\End(TM))$ and $\rho\in C^{\infty}(M)$  satisfying \eqref{hyp1} and \eqref{hyp2} determine a unique conformal infinitesimal bending of $f$. 
\end{theorem}

\begin{remarks}\label{re:infvar} \emph{$1)$ For an infinitesimal variation of a hypersurface  $f\colon M^n\to \mathbb{R}^{n+1}$, its associated tensor $\mathcal{B}$ satisfies \eqref{hyp1} with $\rho=0$, and \eqref{hyp2} reduces to the Codazzi equation for $\mathcal{B}$. \vspace{1ex}\\
$2)$ By Proposition $12$ in \cite{DJV} (respectively, Theorem $13$ in \cite{DV}), a conformal infinitesimal bending (respectively, infinitesimal bending) of a conformal  infinitesimal variation (respectively,  infinitesimal variation) of a hypersurface $f\colon M^n\to \R^{n+1}$, $n\geq 3$, is trivial if and only if its associated tensor $\mathcal{B}$ has the form $\mathcal{B}=\varphi I$ for some $\varphi\in C^{\infty}(M)$ (respectively, its  associated tensor $\mathcal{B}$ vanishes). }
\end{remarks}

\section{Proof of Theorem \ref{prop:cibxmib}}

This section is devoted to the proof of Theorem \ref{prop:cibxmib}, for which we first establish several preliminary facts. \vspace{1ex}

Let $f\colon M^n\to \R^m$ be an isometric immersion and let $F\colon (-\epsilon, \epsilon)\times M^n\to\R^m$ be a smooth variation of $f$ by immersions $f_t=F(t,\cdot)$ with $f_0=f$. From now on, given a  one-parameter family of vector fields  $X^t\in \mathfrak{X}(M)$, we define 
$X'\in \mathfrak{X}(M)$ by setting, for each $x\in M^n$, 
$$
X'(x)=\frac{\partial}{\partial t}|_{t=0}X^t(x).
$$
For the proofs of the next two lemmas we refer to  \cite{JT} (see Lemma~$4$ and Lemma~$5$ therein, respectively).

\begin{lemma}\label{basic} For any fixed $x\in M^n$, the velocity vector at $t=0$ of the smooth curve $t\mapsto f_{t*}X^t(x)$ is
$$
\frac{\partial}{\partial t}|_{t=0}f_{t*} X^t(x)=\tilde \nabla_{X(x)}\T+f_*X'(x),
$$
where $\T$ is the variational vector field of $F$. 
\end{lemma}

\begin{lemma}\label{basic2} If  $\a^t$ denotes the second fundamental form of $f_t$, then 
$$
\<\frac{\partial}{\partial t}|_{t=0}\a^t(X,Y), \eta\>=\<\nab_X\nab_Y\T-\nab_{\nabla_XY}\T,\eta\>
$$
for all $X,Y\in\mathfrak{X}(M)$ and $\eta\in \Gamma(N_fM)$. 
\end{lemma}

Taking into account \eqref{eq:beta}, Lemma \ref{basic2} yields the following for a conformal infinitesimal variation of $f$.

\begin{corollary}\label{basicc3}
If $F$ is a conformal infinitesimal variation of $f$ and $\rho\in C^{\infty}(M)$  is the conformal factor associated to its  conformal infinitesimal bending $\T$, then
\be\label{deranor}
\<\frac{\partial}{\partial t}|_{t=0}\a^t(X,Y),\eta\>=\<\beta(X,Y)+\rho\a(X,Y),\eta\>
\ee
for all $X,Y\in\mathfrak{X}(M)$ and $\eta\in\Gamma(N_fM)$. 
\end{corollary}

For a conformal infinitesimal variation $F$ of $f$ and a (local) orthonormal frame $\{X_i\}_{1\leq i\leq n}$ with respect to the metric induced by $f$, let $X_i^t\in\mathfrak{X}(M)$, $1\leq i\leq n$, $t\in I$, be a smooth one-parameter family of tangent frames such that $X_i^0=X_i$, $1\leq i\leq n$, and
$
\<f_{t*}X_i^t,f_{t*}X_j^t\>=\delta_{ij}
$
for all $1\leq i,j\leq n$ and $t\in I$, that is, $\{X_i^t\}_{1\leq i\leq n}$ is an orthonormal frame for the metric induced by $f_t=F(t,\cdot)$. 

\begin{lemma}
The vector fields $X_i'$, $1\leq i\leq n$, satisfy
\be\label{xilinhai}
\<X_i',X_i\>=-\rho
\ee
and
\be\label{xilinhak}
\<X_i',X_j\>+\<X_i,X_j'\>=0
\ee
for all $1\leq i,j\leq n$ with $i\neq j$.
\end{lemma}
\proof Taking the derivative with respect to $t$ of
$
\<f_{t*}X_i^t,f_{t*}X_j^t\>=\delta_{ij}
$
at $t=0$ and using Lemma \eqref{basic}, we obtain
\begin{align*}
 0= &  \frac{\partial}{\partial t}|_{t=0}\<f_{t*}X_i^t,f_{t*}X_j^t\>\\
     =&\<\nab_{X_i}\T+f_*X_i',f_*X_j\>+\<f_*X_i,\nab_{X_j}\T+f_*X_j'\>\\
    =& \<X_i',X_j\>+\<\nab_{X_i}\T,f_*X_j\>+\<X_i,X_j'\>+\<f_*X_i,\nab_{X_j}\T\>.
\end{align*}

Combining the preceding equation with \eqref{bending} yields
$$
\<X_i',X_j\>+\<X_i,X_j'\>+2\rho\<X_i,X_j\>=0,\,\,\,1\leq i,j\leq n. \qed
$$

\begin{lemma}
Let $f\colon M^n\to\R^m$ be an isometric immersion and let $F\colon (-\epsilon, \epsilon)\times M^n\to\R^m$ be a conformal infinitesimal variation of $f$ with corresponding conformal infinitesimal bending $\T$ and conformal factor $\rho\in C^{\infty}(M)$.
Let $\phi_t$ be given by \eqref{phi} for each immersion $f_t=F(t,\cdot)$, $t\in (-\epsilon, \epsilon)$. Then  
\be\label{derphi}
\frac{\partial}{\partial t}|_{t=0}\phi_t^2=2n\Delta\rho-2\phi^2\rho+2n^2\<\mathcal{L},\mathcal{H}\>,
\ee
where $\mathcal{H}$ is the mean curvature vector field of $f$ and $\mathcal{L}\in\Gamma(N_fM)$ is given by \eqref{meanbeta}.
\end{lemma}
 \proof 
Let us first compute $\partial/\partial t|_{t=0}\|\a^t\|^2$.
Let $\{X_i^t\}_{1\leq i\leq n}$ be a one-parameter family of frames such that, for each fixed $t\in (-\epsilon, \epsilon)$, $\{X_i^t\}_{1\leq i\leq n}$ is orthonormal with respect to the metric induced by $f_t$. By \eqref{deranor} we have
\begin{equation}\label{ddt}
\<\frac{\partial}{\partial t}|_{t=0}\a^t(X^t_i,X^t_j), \eta\>
    =\<\beta(X_i,X_j)+\rho\a(X_i,X_j)+\a(X_i',X_j)+\a(X_i,X_j'),\eta\>,
\end{equation}
hence
\begin{align*}
    \frac{1}{2}\frac{\partial}{\partial t}|_{t=0}\|\a^t(X^t_i,X^t_j)\|^2=&\<\frac{\partial}{\partial t}|_{t=0}\a^t(X^t_i,X^t_j),\a(X_i,X_j)\>\\
    =&\<\beta(X_i,X_j)+\rho\a(X_i,X_j),\a(X_i,X_j)\>\\
    &+\<\a(X_i',X_j)+\a(X_i,X_j'),\a(X_i,X_j)\>.
\end{align*}
Thus
\begin{align*}
    \frac{\partial}{\partial t}|_{t=0}\|\a^t\|^2=&2\rho\|\a\|^2+2\sum_{i,j=1}^n\<\beta(X_i,X_j),\a(X_i,X_j)\> \\
    &+2\sum_{i,j=1}^n\<\a(X_i',X_j)+\a(X_i,X_j'),\a(X_i,X_j)\>\\
=&2\rho\|\a\|^2+2\sum_{i,j=1}^n\<\beta(X_i,X_j),\a(X_i,X_j)\> \\
&+4\sum_{i,j=1}^n\<\a(X_i',X_j),\a(X_i,X_j)\>.
\end{align*}
It follows from \eqref{derGaussC} that
 \begin{align*}
     2\<\beta(X_i,X_j),\a(X_i,X_j)\>=&
\<\beta(X_i,X_i),\a(X_j,X_j)\>+\<\beta(X_j,X_j),\a(X_i,X_i)\>\\
&+\<X_j,X_j\>\hess\rho(X_i,X_i)+\<X_i,X_i\>\hess\rho(X_j,X_j)
\end{align*}
for all $1\leq i, j\leq n$ with $i\neq j$. 
 Therefore
 $$
     2\sum_{i,j=1}^n\<\beta(X_i,X_j),\a(X_i,X_j)\>= 2(n-1)\Delta\rho
     +2\sum_{i, j=1}^n\<\beta(X_i,X_i),\a(X_j,X_j)\>,
 $$
 where $\Delta$ denotes the Laplacian.
 On the other hand, from the Gauss equation,
 $$
     \<\a(X_i',X_j),\a(X_i,X_j)\>=\<\a(X_i',X_i),\a(X_j,X_j)\>+\<R(X_i',X_j)X_i,X_j\>,
 $$
 where $R$ denotes the Riemann curvature tensor of $M^n$, we obtain 
 $$
     \sum_{i\neq j}\<\a(X_i',X_j),\a(X_i,X_j)\>=
     \sum_{i\neq j}\<\a(X_i',X_i),\a(X_j,X_j)\>-\sum_{i=1}^n \Ric(X_i',X_i).
 $$
 Thus
 \begin{align*}
     \frac{\partial}{\partial t}|_{t=0}\|\a^t\|^2=&2\rho\|\a\|^2+2(n-1)\Delta\rho
     +2\sum_{i, j=1}^n\<\beta(X_i,X_i),\a(X_j,X_j)\>\\
     &+4\sum_{i,j=1}^n\<\a(X_i',X_i),\a(X_j,X_j)\>-4\sum_{i=1}^n \Ric(X_i',X_i)\\
=&2\rho\|\a\|^2+2(n-1)\Delta\rho+2n^2\<\mathcal{L},\mathcal{H}\>\\
&+4n\sum_{i=1}^n\<\a(X_i',X_i),\mathcal{H}\>-4\sum_{i=1}^n\Ric(X_i',X_i).
 \end{align*}
By \eqref{xilinhai}, we may write $X_i'=-\rho X_i+\sum_{i\neq k}\<X_i',X_k\>X_k$; hence
\begin{align*}
    \sum_{i=1}^n\<\a(X_i',X_i),\mathcal{H}\>=&-\rho n\|\mathcal{H}\|^2+\sum_{i\neq k}\<X_i',X_k\>\<\a(X_k,X_i),\mathcal{H}\>\\
    =&-\rho n\|\mathcal{H}\|^2,
\end{align*}
where the last equality follows from \eqref{xilinhak}.
Similarly,
$$
    \sum_{i=1}^n\Ric(X_i',X_i)=-\rho n(n-1)s,
$$
where 
$
{\displaystyle s=\frac{1}{n(n-1)}\sum_{i=1}^n\Ric(X_i,X_i)}
$
 is the scalar curvature of $M^n$.
Thus 
\begin{align*}
     \frac{\partial}{\partial t}|_{t=0}\|\a^t\|^2=&2\rho\|\a\|^2+2(n-1)\Delta\rho+2n^2\<\mathcal{L},\mathcal{H}\>\\
     &-4n^2\rho\|\mathcal{H}\|^2+4\rho n(n-1)s.
\end{align*}
Using that
$$
s=\frac{n}{n-1}\|\mathcal{H}\|^2-\frac{1}{n(n-1)}\|\a\|^2,
$$
we obtain
\be\label{um}
\frac{\partial}{\partial t}|_{t=0}\|\a^t\|^2=2(n-1)\Delta \rho+2n^2\<\mathcal{L},\mathcal{H}\>-2\rho\|\a\|^2.
\ee
We now compute $\partial/\partial t|_{t=0}\|\mathcal{H}^t\|^2$. With $\{X_i^t\}_{1\leq i\leq n}$ as above, we have
\begin{align*}
    \frac{\partial}{\partial t}|_{t=0}\|\mathcal{H}^t\|^2=&2\<\frac{\partial}{\partial t}|_{t=0}\mathcal{H}^t,\mathcal{H}\>\\
    =&2\<\frac{1}{n}\sum_{i=1}^n\frac{\partial}{\partial t}|_{t=0}\a^t(X^t_i,X^t_i),\mathcal{H}\>\\
    =&2\rho\|\mathcal{H}\|^2+2\<\mathcal{L},\mathcal{H}\>+\frac{4}{n}\sum_{i=1}^n\<\a(X_i',X_i),\mathcal{H}\>,
\end{align*}
where the last step follows from \eqref{ddt}.
Using  \eqref{xilinhai} and \eqref{xilinhak} as before, we obtain
\be\label{dois}
    \frac{\partial}{\partial t}|_{t=0}\|\mathcal{H}^t\|^2=2\<\mathcal{L},\mathcal{H}\>-2\rho\|\mathcal{H}\|^2.
\ee
It follows from \eqref{um} and \eqref{dois} that
\begin{align*}
    \frac{\partial}{\partial t}|_{t=0}\phi_t^2=&\frac{n}{n-1}\left(2(n-1)\Delta \rho+2n^2\<\mathcal{L},\mathcal{H}\>-2\rho\|\a\|^2\right.\\
    &\left.-2n\<\mathcal{L},\mathcal{H}\>+2n\rho\|\mathcal{H}\|^2\right)\\
    =&\frac{n}{n-1}(2(n-1)\Delta\rho+2n(n-1)\<\mathcal{L},\mathcal{H}\>\\
    &-2\rho\|\a\|^2+2\rho n\|\mathcal{H}\|^2)\\
    =&2n\Delta\rho+2n^2\<\mathcal{L},\mathcal{H}\>-2\phi^2\rho,
\end{align*}
where we have used \eqref{phi} in the last equality. \vspace{2ex}
\qed

\noindent \emph{Proof of Theorem \eqref{prop:cibxmib}:}
If $\T$ is the variational vector field of an  infinitesimal Moebius variation $F\colon (-\epsilon, \epsilon)\times M^n\to \mathbb{R}^m$ of $f$, then the corresponding conformal factor $\rho$ is given by \eqref{rhomoeb}. 
Thus  \eqref{derphi} yields
$$
-2\phi^2\rho=2n\Delta\rho-2\phi^2\rho+2n^2\<\mathcal{L},\mathcal{H}\>,
$$
and hence \eqref{eqrhomoeb} holds.

For the converse, assume that $\T$ is a conformal infinitesimal bending of $f$ whose conformal factor $\rho$ satisfies \eqref{eqrhomoeb}.
The variation $\F\colon\R\times M^n\to\R^m$ given by
$
\F(t,x)=f(x)+t\T(x)
$
is a conformal infinitesimal variation with variational
vector field $\T$.  Let $f_t=\F(t,\cdot)$ and let $\phi_t$ be given by \eqref{phi} for each $f_t$, $t\in \R$. We claim that $\F$ is an   infinitesimal Moebius variation of $f$. Indeed, we have  
$$
    \frac{\partial}{\partial t}|_{t=0}(\phi_t^2\<f_{t*}X,f_{t*}Y\>)=
    \frac{\partial}{\partial t}|_{t=0}(\phi_t^2)\<X,Y\>+\phi^2(\<\tilde \nabla_X\T,f_*Y\>+\<f_*X,\tilde\nabla_Y\T\>)
$$
for all $X,Y\in\mathfrak{X}(M)$, hence 
$$
    \frac{\partial}{\partial t}|_{t=0}(\phi_t^2\<f_{t*X},f_{t*}Y\>)=\left(\frac{\partial}{\partial t}|_{t=0}(\phi_t^2)+2\phi^2\rho\right)\<X,Y\>
    $$
by \eqref{bending}.  On the other hand, from \eqref{derphi} and \eqref{eqrhomoeb} we have
$$
    \frac{\partial}{\partial t}|_{t=0}(\phi_t^2)+2\phi^2\rho=0,
$$
which  proves the claim and completes the proof. 
\vspace{2ex}\qed

Before concluding this section, we state for later use the following consequence of some of the preceding computations (see  \cite{DV} for the corresponding fact for (isometric) infinitesimal variations).
 
 \begin{proposition}\label{eq:alinha} Let $F\colon (-\epsilon, \epsilon)\times M^n\to\R^{n+1}$ be a conformal infinitesimal variation of an isometric immersion $f\colon M^n\to \R^{n+1}$. Let $N_t$ be a unit vector field normal to  $f_t=F(t,\cdot)$, $t\in (-\epsilon, \epsilon)$, and denote by $A_t$  the corresponding shape operator. Then the tensor $\mathcal{B}\in \Gamma(\End(TM))$ associated with $F$ satisfies
 \be\label{BAlinha}
\mathcal{B}=A'+\rho A,
\ee
where $\rho\in C^{\infty}(M)$ is the conformal factor of $F$ and 
$A'=\partial/\partial t|_{t=0}A_t$.
\end{proposition}
\proof 
It follows from \eqref{bending} and Lemma \ref{basic} that
\begin{eqnarray*}
\partial/\partial t|_{t=0}\<\a^t(X,Y),N_t\>&=&\partial/\partial t|_{t=0}\<{f_t}_*A_tX,
{f_t}_*Y\>\vspace{1ex}\\
&=&\<A'X,Y\>+2\rho\<AX,Y\>
\end{eqnarray*}
for all $X,Y\in\mathfrak{X}(M)$.  On the other hand, from 
\eqref{deranor} we obtain
$$
\partial/\partial t|_{t=0}\<\a^t(X,Y),N_t\>=\<\mathcal{B}X,Y\>+\rho\<AX,Y\>.
$$
Comparing the two preceding equations yields \eqref{BAlinha}. \qed

\section{Proof of Theorem \ref{thm:main}}

In this section we prove Theorem \ref{thm:main}. We start with some  preliminary results, 
 which  make use of the following lemma in \cite{DJV} (see Lemma~$14$ therein).

\begin{lemma}\label{djv} Let  $f\colon M^n\to\R^{n+1}$, $n\geq 5$, be a hypersurface that admits a conformal infinitesimal  bending $\T$ that is non-trivial on any open subset. Then its associated tensor $\mathcal{B}$, the Hessian $H$ of its conformal factor $\rho$, and the shape operator $A$ of $f$ share, on an open and dense subset of $M^n$, a common eigenbundle $\Delta$ of constant dimension $\dim \Delta\geq n-2$. 
\end{lemma}  

 The next result states how Theorem  \eqref{prop:cibxmib} reads for hypersurfaces.

\begin{proposition} \label{prop:imbhyp}
Let  $f\colon M^n\to\R^{n+1}$, $n\geq 5$, be an umbilic-free infinitesimally  Moebius bendable hypersurface.   Then there exists an open and dense subset $\mathcal{U}$ of $M^n$ where $\mathcal{B}$, $H$ and $A$ share a common eigenbundle $\Delta$ of rank $n-2$ and  $\trace(A-\lambda I)\trace(\mathcal{B}-bI)=0$, 
where $b, \lambda \in C^{\infty}(\mathcal{U})$ are such that $\mathcal{B}|_\Delta=bI$ and $A|_\Delta=\lambda I$.
\end{proposition}
\proof  Since $f$ is  infinitesimally  Moebius bendable, it admits, in particular, a conformal infinitesimal  bending $\T$ that is non-trivial on any open subset.  By Lemma 
\ref{djv}, there exists an open and dense subset $\mathcal{U}$ of $M^n$ where $\mathcal{B}$, $H$ and $A$ share a common eigenbundle $\Delta$ of constant dimension $\dim \Delta\geq n-2$. 
In the proof of Proposition 15 in \cite{DJV} (see Eq. $(35)$ therein), it was shown that
\be\label{suma}
bA+\lambda(\mathcal{B}-bI)+\hess\rho=0
\ee
on $\mathcal{U}$. Let $\mathcal{H}$ and $\mathcal{L}$ be given by
$\trace A=n\mathcal{H}$ and $\trace \mathcal{B}=n\mathcal{L}$. 
Taking traces in \eqref{suma} yields
$$
nb\mathcal{H}+n\lambda\mathcal{L}-n\lambda b+\Delta\rho=0.
$$
We write the preceding equation as 
$$
n(\mathcal{H}-\lambda)(\mathcal{L}-b)=\Delta\rho+n\mathcal{L}\mathcal{H},
$$
which is also equivalent to 
\begin{equation}\label{eq:le}
\trace(A-\lambda I)\trace(\mathcal{B}-bI)=n(\Delta\rho+n\mathcal{L}\mathcal{H}).
\end{equation}
Taking into account that \eqref{eqrhomoeb} reduces to $\Delta\rho+n\mathcal{L}\mathcal{H}=0$,
it follows from \eqref{eq:le} and Theorem~\eqref{prop:cibxmib} that $\trace(A-\lambda I)\trace(\mathcal{B}-bI)~=~0$. 

Finally, notice that the preceding condition can not occur on any open subset  where $\dim \Delta= n-1$. Indeed, if  $\dim \Delta= n-1$, then the condition   
 $\trace(A-\lambda I)=0$ would imply that $A=\lambda I$, whereas $\trace(\mathcal{B}-bI)=0$
 would yield $\mathcal{B}=bI$, in contradiction with the assumptions that $f$ is free of umbilic points and that the infinitesimal bending $\T$ is non-trivial, respectively. \qed
 
 \begin{lemma}\label{lem:umb}
The distribution $\Delta$ given by Proposition \ref{prop:imbhyp} is umbilical.
\end{lemma}

\proof If $\Delta=\ker (A-\lambda I)$, then it is the eigenbundle corresponding to the principal curvature $\lambda$, and hence umbilical. Thus we only need to consider the case in which $\Delta$ coincides with $\ker (B-bI)$ and is a proper subspace of $\ker(A-\lambda I)$.
Equation \eqref{hyp2} can be written as
$$
(\nabla_X(\mathcal{B}-bI))Y-(\nabla_Y(\mathcal{B}-bI))X+(X\wedge Y)(A\nabla\rho-\nabla b)=0
$$
for any $X,Y\in\mathfrak{X}(M)$, where $\nabla b$ and $\nabla\rho$ are the gradients of $b$ and of the conformal factor $\rho$, respectively.
Let $T,S\in\Gamma(\Delta)$ be orthogonal and take $X\in\Gamma(\Delta^\perp)$. Evaluating  the preceding equation in $X$ and $T$ and taking the inner product of both sides with $S$ gives
$
\<\nabla_T(\mathcal{B}-bI)X,S\>=0.
$
Since we are assuming that $\rank(\mathcal{B}-bI)=2$,  the above equation gives 
$$
(\nabla_TS)_{\Delta^\perp}=0
$$
for all $T,S\in\Gamma(\Delta)$ with $\<T,S\>=0$. Thus $\Delta$ is an umbilical distribution.
\vspace{1ex}\qed

 From now on, for  $\mathcal{U}$ and $b, \lambda \in C^{\infty}(\mathcal{U})$ as in Proposition \ref{prop:imbhyp}, we denote 
$\bar{A}=A-\lambda I$ and $\bar{\mathcal{B}}=\mathcal{B}-bI$.
 
 \begin{proposition} \label{prop2:imbhyp}
Let  $f\colon M^n\to\R^{n+1}$, $n\geq 5$, be an umbilic-free infinitesimally  Moebius bendable hypersurface,  let $\mathcal{U}$ be the open and dense subset of $M^n$ given by Proposition \ref{prop:imbhyp}, and let   $\mathcal{U}_1$ be the subset of $\mathcal{U}$ where $\trace \bar{\mathcal{B}}=0$. Then $\mathcal{U}_1=\mathcal{Y}_1\cup \mathcal{Y}_2$, where the following holds on $\mathcal{Y}_1$ and $\mathcal{Y}_2$, respectively:
\begin{itemize}
\item[(i)] $\bar{A}|_{\Delta^\perp}$ is a multiple of the identity endomorphism $I\in \Gamma(\End(\Delta^\perp))$; 
 \item[(ii)] there exists at each point an orthonormal basis $\{X,Y\}$ of $\Delta^{\perp}$ given by principal directions of $f$ and $\theta\in \mathbb{R}$ such that $\bar{\mathcal{B}}X=\theta Y$ and $\bar{\mathcal{B}}Y=\theta X$.
 \end{itemize}
\end{proposition}
\proof  It follows from   \eqref{suma}  that, at each $x\in \mathcal{U}$, Eq. \eqref{hyp1} can be written as
\begin{equation}\label{hyp1b}
\bar{\mathcal{B}}X \wedge \bar A Y= \bar{\mathcal{B}}Y \wedge \bar A X,
\end{equation}
or equivalently,
$$\<\bar{A}Y,X\>\<\bar{\mathcal{B}}X, Y\>-\<\bar{\mathcal{B}}X,X\>\<\bar{A}Y, Y\>=\<\bar{A}X,X\>\<\bar{\mathcal{B}}Y,Y\>-\<\bar{\mathcal{B}}Y,X\>\<\bar{A}X,Y\>$$
for all $X, Y\in T_x \mathcal{U}$. 
Applying the preceding equation to orthogonal unit eigenvectors $X$ and $Y$ of $\bar A|_{\Delta^\perp}$, with $\bar A X=\mu_1X$ and $\bar A Y=\mu_2Y$, gives 
$$
-\mu_2\<\bar{\mathcal{B}}X,X\>=\mu_1\<\bar{\mathcal{B}}Y,Y\>.
$$

 Therefore, at each point of $\mathcal{U}_1$, either $\mu_1=\mu_2:=\mu\neq 0$, and hence $\bar A|_{\Delta^\perp}=\mu I$, or $\<\bar{\mathcal{B}}X,X\>=0=\<\bar{\mathcal{B}}Y,Y\>$.  In the latter case, denoting
$\theta=\<\bar{\mathcal{B}}X,Y\>$, we have $\bar{\mathcal{B}}X=\theta Y$ and $\bar{\mathcal{B}}Y=\theta X$. \vspace{2ex}\qed

Given a distribution $\Delta$ on a Riemannian manifold $M^n$, recall that the 
\emph{splitting tensor}
$C\colon\Gamma(\Delta)\to\Gamma(\End(\Delta^\perp))$ of $\Delta$
is defined by 
$$
C_TX=-\nabla_X^hT
$$ 
for all $T\in\Gamma(\Delta)$ and $X\in\Gamma(\Delta^\perp)$, where $\nabla_X^hT=(\nabla_XT)_{\Delta^\perp}$.

\begin{proposition}\label{prop:class}
Let  $f\colon M^n\to\R^{n+1}$, $n\geq 5$, be an umbilic-free infinitesimally  Moebius bendable hypersurface carrying a principal curvature of constant multiplicity $n-2$ with 
corresponding eigenbundle $\Delta$. Assume that at no point of $M^n$ the splitting tensor $C\colon\Gamma(\Delta)\to\Gamma(\End(\Delta^\perp))$  of $\Delta$ satisfies $C(\Gamma(\Delta))\subset \spa\{I\}$. Then the central sphere congruence of $f$  is 
determined by a minimal space-like surface $s\colon L^2\to \Sf_1^{n+2}$. 
\end{proposition}
\proof Denoting by $\lambda$ the principal curvature of $f$ with constant multiplicity $n-2$ with respect to a unit normal vector field $N$, the map $$x\in M^n\mapsto f(x)+\frac{1}{\lambda(x)}N$$ determines a two-parameter congruence of hyperspheres that is enveloped by $f$. As explained in the introduction, this congruence of hyperspheres is determined by a space-like surface $s\colon L^2\to\Sf_1^{n+2}$.

Since $f$ is  infinitesimally  Moebius bendable, it admits, in particular, a conformal infinitesimal  bending that is non-trivial on any open subset. By Proposition $15$ in \cite{DJV}, the hypersurface $f$ is either elliptic, hyperbolic or parabolic with respect to $J\in\Gamma(\End(\Delta^\perp))$ satisfying  $J^2=-I$, $J^2=I$ or $J^2=0$, respectively, with $J\neq I$ if $J^2=I$ and  $J\neq 0$ if $J^2=0$. Moreover, the tensor $\mathcal{B}$ associated with $\T$ satisfies
\be\label{B-b}
\bar{\mathcal{B}}=\mu \bar{A} J,
\ee
where $0\neq\mu\in C^\infty(M)$ is constant along the leaves of $\Delta$.

It was also shown in \cite{DJV} that, in the hyperbolic and elliptic cases, the tensor $J$ is projectable with respect to the quotient map   $\pi\colon M^n\to L^2$ onto the spaces of leaves of the eigenbundle $\Delta$ of $\lambda$, that is, there exists 
  $\bar{J} \in \mbox{End}(TL)$ such that $\bar{J}\circ \pi_*=\pi_*\circ J$. Moreover,  the surface $s\colon L^2\to\Sf_1^{n+2}$ is either a  \emph{special elliptic} or  \emph{special hyperbolic} surface with respect to $\bar J$. This means that 
\be\label{sffsur}
\alpha^s(\bar{J}\bar{X},\bar{Y})=\alpha^s(\bar{X},\bar{J}\bar{Y})
\ee
for all $\bar{X},\bar{Y}\in \mathfrak{X}(L)$, and that there exists $\mu\in C^\infty(L)$ such that $\mu \bar J$ is a Codazzi tensor on $L^2$.

  In the sequel we will show that, under the assumptions of the proposition, the tensors $J$  and $\bar{J}$ act as a rotation of angle $\pi/2$ on $\Delta^\perp$ and on each tangent space of $L^2$, respectively, that is, both $J$  and $\bar{J}$ are orthogonal tensors satisfying 
  $J^2=-I$ and $\bar{J}^2=-I$. From the orthogonality of $J$ and the symmetry of $\bar{B}$ it will follow that the tensor $\bar{A}=A-\lambda I$ is traceless by \eqref{B-b}. This implies that $\lambda$ is the mean curvature function of $f$, and hence the congruence of hyperspheres determined by $s$ is its central sphere congruence. On the other hand,  the orthogonality of $\bar{J}$ and the fact that $\bar{J}^2=-I$ implies the minimality of $s$ by 
\eqref{sffsur}, and this will conclude the proof. 

   First we rule out the parabolic case. So, assume that there exists $J\in\Gamma(\Delta^\perp)$ such that $J^2=0$, $J\neq 0$, $\nabla^h_T J=0$ for all $T\in\Gamma(\Delta)$, and such that  $C_T\in\spa\{I,J\}$ for all $T\in\Gamma(\Delta)$. By Proposition $16$ of \cite{DJV}, 
 $f$ is  conformally ruled, with the leaves of the distribution $\Delta\oplus\ker(J)$ as the rulings of $f$.  
Let $X,Y\in\Gamma(\Delta^\perp)$ be an orthonormal basis of $\Delta^\perp$ such that $JX=Y$ and $JY=0$, let $\lambda_1, \lambda'\in C^{\infty}(M)$ be such that $\bar A X=\lambda_1X+\lambda' Y$ and $\bar A  Y=\lambda' X$.
From \eqref{B-b} we see that $\bar{\mathcal B}X=\mu\lambda' X$ and $\bar{\mathcal B} Y=0$.
Since  $\mathcal{B}$ is not a multiple of the identity endomorphism, then $\lambda'\neq 0$, and hence $\trace \bar{\mathcal B}\neq 0$. 
It follows from Proposition \eqref{prop:imbhyp} that $\lambda_1=\trace \bar A=0$.
It follows from the Codazzi equation  that
\be\label{codA}
\nabla^h_T \bar A=\bar{A}C_T
\ee
for any $T\in\Gamma(\Delta)$.
For a fixed $T\in\Gamma(\Delta)$, write
$
C_T=dI+eJ
$ 
for some smooth functions $d$ and $e$.
On one hand, $
\nabla^h_T\bar{A}X=\nabla^h_T\lambda'Y=T(\lambda')Y,
$
where we have used that $\nabla^h_TY=0$, for $Y$ is tangent to the rulings.
On the other hand, 
$$ \bar{A}C_TX=\bar{A}(dX+eY)=d\lambda'Y+e\lambda'X.$$
Therefore  $e=0$ by \eqref{codA} and, since $T\in\Gamma(\Delta)$ was chosen arbitrarily, it follows that $C_T\in \spa\{I\}$ for any $T\in\Gamma(\Delta)$, a contradiction with our assumption.

Now assume that $f$ is hyperbolic, that is, that there exists 
$J\in\Gamma(\End(\Delta^\perp))$ such that  $J^2=I$, with $J\neq I$, $\nabla^h_T J=0$ and 
such that $C_T\in\spa\{I,J\}$ for all $T\in\Gamma(\Delta)$.
Let $\{X,Y\}$ be a frame of $\Delta^\perp$ of unit eigenvectors of $J$, with $JX=X$ and $JY=-Y$.
Since $\nabla^h_T J=0$  for all $T\in\Gamma(\Delta)$, it follows that 
 $\nabla^h_T X=0=\nabla^h_T Y$. The symmetry of $\bar{\mathcal{B}}=\mu\bar{A}J$ yields
 $\<\bar{A}X,Y\>=0$. Write
 $\bar{A}X=\alpha X+\beta Y$ and $\bar{A}Y=\gamma X+ \delta Y$
 for some smooth functions $\alpha$, $\beta$,  $\gamma$,  $\delta$. Then
\begin{equation}\label{eq:bg}
\<\bar{A}X,X\>=\alpha+\beta\<Y,X\>,\,\,\,\,\,\<\bar{A}Y,Y\>=\gamma\<Y,X\>+\delta,\end{equation}
and from $\<\bar{A}X,Y\>=0=\<X, \bar{A}Y\>$ we obtain
 \begin{equation}\label{eq:ad}
 \alpha\<X,Y\>+\beta=0=\gamma+\delta\<X,Y\>.
 \end{equation}
 On the other hand,  writing as before $C_T=dI+eJ$
for some smooth functions $d$ and $e$,
 Eq. \eqref{codA} gives
\begin{equation}\label{eq:n1}
 \<\nabla^h_T\bar{A}X, X\>=\<\bar{A}C_TX, X\>=(d+e)\<\bar{A}X, X\>,
 \end{equation}
 and similarly,
 \begin{equation}\label{eq:n2}
 \<\nabla^h_T\bar{A}Y, Y\>=\<\bar{A}C_TY, Y\>=(d-e)\<\bar{A}Y, Y\>.
 \end{equation}

 Suppose that $\trace\bar{A}=0$. Then  $\alpha=-\delta$, hence \eqref{eq:ad}  implies that $\beta=-\gamma$. Thus
 $\<\bar{A}X,X\>=-\<\bar{A}Y,Y\>$
 by \eqref{eq:bg}, and hence 
 $$\<\nabla^h_T\bar{A}X, X\>=T\<\bar{A}X, X\>=-T\<\bar{A}Y, Y\>=-\<\nabla^h_T\bar{A}Y, Y\>.
 $$
 Comparing with \eqref{eq:n1} and  \eqref{eq:n2} gives $d+e=d-e$,  for $\<\bar{A}X, X\>\neq 0$ by the assumption that $\rank\bar A=2$.
  Hence $e=0$. 
 
 If $\trace \bar{\mathcal{B}}=0$, then \eqref{B-b} gives $\alpha=\delta$, and hence $\gamma=\beta$ by \eqref{eq:ad}. Therefore 
 $\<\bar{A}X,X\>=\<\bar{A}Y,Y\>$
 by \eqref{eq:bg}, and hence 
 $$\<\nabla^h_T\bar{A}X, X\>=T\<\bar{A}X, X\>=T\<\bar{A}Y, Y\>=\<\nabla^h_T\bar{A}Y, Y\>.
 $$
 Then we obtain as before that $e=0$ by comparing  with \eqref{eq:n1} and  \eqref{eq:n2}, and we conclude as in the parabolic case that $C_T\in \spa\{I\}$ for any $T\in\Gamma(\Delta)$, a contradiction.

Finally, suppose that $f$ is elliptic, that is, that there exists 
$J\in\Gamma(\End(\Delta^\perp))$ such that  $J^2=-I$, $\nabla^h_T J=0$, and such that
$C_T\in\spa\{I,J\}$ for all $T\in\Gamma(\Delta)$.
Let $\{X,Y\}$ be a frame of $\Delta^\perp$  such that $JX=Y$ and $JY=-X$.
This is equivalent to asking the complex vector fields $X-iY$ and $X+iY$ to be pointwise  eigenvectors of the $\mathbb{C}$-linear extension of $J$,  also denoted by $J$,  associated to the eigenvalues $i$ and $-i$, respectively.  Thus 
$z(X-iY)=(sX+tY)+i(tX-sY)$ and $z(X+iY)=(sX-tY)+i(tX+sY)$ are also eigenvectors of $J$  associated to $i$ and $-i$, respectively, for any $z=s+it\in \mathbb{C}$, that is,   $\bar{X}=sX+tY$ and $\bar{Y}=-tX+sY$ form a new frame of $\Delta^\perp$  such that $J\bar X=\bar Y$ and $J\bar Y=-\bar X$. It is easily seen that $s$ and $t$ can be chosen so that $\bar{X}$ and $\bar{Y}$ are unit vector fields. In summary, we can always choose a frame $\{X,Y\}$ of \emph{unit} vector fields such that $JX=Y$ and $JY=-X$.

Since $\nabla^h_T J=0$  for all $T\in\Gamma(\Delta)$, then $J\nabla^h_T X=\nabla^h_T Y$ and $J\nabla^h_T Y=-\nabla^h_T X$. Denoting $\hat X=\nabla^h_T X$ and $\hat Y=\nabla^h_T Y$, it follows that
 $\hat X-i\hat Y=(s+it)(X-iY)$
 for some $s+it\in \mathbb{C}$, that is,
 $$\hat X=sX+tY\,\,\,\,\,\mbox{and}\,\,\,\,\,\hat Y=-tX+sY.$$ 
 Since $X$ and $Y$ have unit length, then $\<\hat X,X\>=0=\<\hat Y, Y\>$. Thus
 $$s+t\<X,Y\>=0=s-t\<X,Y\>,$$
 and hence $t\<X,Y\>=0=s$. 
 
 Assume that $J$ is not an orthogonal tensor, that is, that $\<X,Y\>\neq 0$. Then 
 $s=0=t$, that is, $\nabla^h_T X=0=\nabla^h_T Y$ for all $T\in \Gamma(\Delta)$.
 
  Write $\bar{A}X=\alpha X+\beta Y$ and $\bar{A}Y=\gamma X+ \delta Y$
 for some smooth functions $\alpha$, $\beta$,  $\gamma$ and $\delta$.
 The symmetry of $\bar{\mathcal{B}}$ gives
 \be\label{eq:sym}
 \<\bar{A}X,X\>+\<\bar{A}Y,Y\>=0.
 \ee
 Then
\begin{equation}\label{eq:bg2}
\<\bar{A}X,X\>=\alpha+\beta\<Y,X\>,\,\,\,\,\,\<\bar{A}Y,Y\>=\gamma\<Y,X\>+\delta,\end{equation}
and from \eqref{eq:sym} and $\<\bar{A}X,Y\>=\<X, \bar{A}Y\>$ we obtain, respectively,
\begin{equation}\label{eq:sym2}
 (\alpha+\delta) +(\beta+\gamma)\<X,Y\>=0
 \end{equation}
 and
 \begin{equation}\label{eq:ad2}
 \alpha\<X,Y\>+\beta=\gamma+\delta\<X,Y\>.
 \end{equation}
 On the other hand, writing as before $
C_T=dI+eJ
$
for some smooth functions $d$ and $e$,
 Eq. \eqref{codA} gives
\begin{equation}\label{eq:n1b}
 \<\nabla^h_T\bar{A}X, Y\>=\<\bar{A}C_TX, Y\>=d\<\bar{A}X, Y\> +e\<\bar{A}Y, Y\>.
 \end{equation}

 Now assume that $\trace\bar{A}=0$. Then  $\alpha=-\delta$, hence $\beta=-\gamma$ by \eqref{eq:sym2}. Thus $\<\bar{A}X,Y\>=0$  by \eqref{eq:ad2}, and hence 
 $\<\nabla^h_T\bar{A}X, Y\>=T\<\bar{A}X, Y\>=0$.
 It follows from \eqref{eq:n1b} that $e=0$,  for $\<\bar{A}Y, Y\>\neq 0$  by the assumption that $\rank\bar A=2$. 
 
 If $\trace\bar{\mathcal{B}}=0$, then \eqref{B-b} gives $\gamma=\beta$, hence $\alpha=\delta$ by \eqref{eq:ad2}. Therefore 
 $\<\bar{A}X,X\>=0=\<\bar{A}Y,Y\>$
 by \eqref{eq:sym2} and \eqref{eq:bg2}.  Then
 $\<\nabla^h_T\bar{A}X, X\>=T\<\bar{A}X, X\>=0$,
 and, on the other hand,
 \begin{eqnarray*}
 \<\nabla^h_T\bar{A}X, X\>&=&\<\bar{A}C_TX,X\>\\
 &=&\<\bar{A}(dX+eY),X\>\\
 &=& d\<\bar{A}X,X\>+e\<\bar{A}Y,X\>\\
 &=& e\<\bar{A}Y,X\>.
 \end{eqnarray*}
 It follows that $e=0$, for $\<\bar{A}Y, X\>\neq 0$ by the assumption that $f$ is free of points with a principal curvature of multiplicity at least $n-1$. 
 We conclude as in the previous cases that $C_T\in \spa\{I\}$ for any $T\in \Gamma(\Delta)$, a contradiction. 

It follows that $J$ must be an orthogonal tensor, that is,  $\<X,Y\>=0$.  It remains to show that the tensor  $\bar{J} \in \mbox{End}(TL)$ given by $\bar{J}\circ \pi_*=\pi_*\circ J$ is also orthogonal. For this, we use the fact that the metric $\<\cdot,\cdot\>'$ on $L^2$ induced by $s$ is related to the metric of $M^n$ by
\be\label{relmetr}
\<\bar{Z},\bar{W}\>'=\<\bar{A}Z, \bar{A}W\>
\ee
for all $\bar{Z},\bar{W}\in\mathfrak{X}(L)$, where $Z$, $W$ are the horizontal lifts of $\bar{Z}$ and $\bar{W}$, respectively. Let $\bar{X}\in\mathfrak{X}(L)$ and denote by $X\in\Gamma(\Delta^\perp)$ its horizontal lift. Using the symmetry of $\bar{A}J$, we have
\begin{align*}
    \<\bar{X},\bar{J}\bar{X}\>'&=\<\bar{A}X, \bar{A}JX\>\\
    &=\<\bar{A}J\bar{A}X,X\>\\
    &=\<J\bar{A}X,\bar{A}X\>\\
    &=0,
\end{align*}
where in the last step we have used that $J$ acts as a rotation of angle $\pi/2$ on $\Delta^\perp$. Using again the symmetry of $\bar{A}J$, the proof of the orthogonality of $\bar{J}$ is completed by noticing that
\begin{align*}
    \<\bar{J}\bar{X},\bar{J}\bar{X}\>'&=\<\bar{A}JX,\bar{A}JX\>\\
    &=\<J\bar{A}JX,\bar{A}X\>\\
    &=\<JJ^t\bar{A}X,\bar{A}X\>\\
    &=-\<J^2\bar{A}X,\bar{A}X\>\\
    &=\<\bar{X},\bar{X}\>'. \qed
\end{align*}    

For the proof of  Theorem \ref{thm:main} 
 we will also need the following fact
(see Theorem $1$ in \cite{DFT} or Corollary $9.33$ in \cite{DT}).

\begin{lemma} \label{le:split} Let $f\colon M^n\to\R^{n+1}$, $n\geq 3$, be a hypersurface
and let $\Delta$ be an umbilical subbundle of rank $n-2$ of the eigenbundle of $f$ correspondent to a principal curvature of $f$.
Then $f$ is  conformally surface-like (with respect to the decomposition $TM=\Delta^\perp\oplus \Delta$) if and only if the splitting tensor $C\colon\Gamma(\Delta)\to\Gamma(\End(\Delta^\perp))$ of $\Delta$ satisfies $C(\Gamma(\Delta))\subset \spa\{I\}$. 
\end{lemma}

\noindent \emph{Proof of Theorem \ref{thm:main}}:  Let $\mathcal{U}$  and $\Delta$ be, respectively,  the open and dense subset of $M^n$ and the distribution of rank $n-2$ given by Proposition \ref{prop:imbhyp}. By that result,  $\mathcal{U}$ splits as $\mathcal{U}=\mathcal{U}_1\cup \mathcal{U}_2$, with  $\trace \bar{\mathcal{B}}=0$ on $\mathcal{U}_1$ and $\trace \bar{A}=0$ on $\mathcal{U}_2$. 

We also consider the decompositions $\mathcal{U}=\mathcal{V}_1\cup \mathcal{V}_2$ and $\mathcal{U}=\mathcal{W}_1\cup \mathcal{W}_2$, where $\mathcal{V}_1$ and $\mathcal{V}_2$ are the subsets  where the dimension of $\ker \bar{A}$ is either $n-2$ or $n-1$, respectively,  $\mathcal{W}_2$ is the 
 subset where $C(\Gamma(\Delta))\subset \spa\{I\}$ and $\mathcal{W}_1=\mathcal{U}\setminus \mathcal{W}_2$.

 In the following we denote by $S^0$ the interior of the subset $S$. We will show that the direct statement holds on the open and dense subset 
 $$\mathcal{U}^*=\mathcal{V}^0_2\cup (\mathcal{V}_1\cap \mathcal{W}_1)\cup  (\mathcal{V}_1\cap \mathcal{W}^0_2\cap \mathcal{U}^0_2)\cup  (\mathcal{V}_1\cap \mathcal{W}^0_2\cap \mathcal{Y}^0_2)\cup (\mathcal{V}_1\cap \mathcal{W}^0_2\cap \mathcal{Y}^0_1),$$
where $\mathcal{Y}_1$ and $\mathcal{Y}_2$ are the subsets of $\mathcal{U}_1$ given by Proposition \ref{prop2:imbhyp}.

It follows from Proposition \ref{prop:class} that the central sphere congruence of $f|_{\mathcal{V}_1\cap \mathcal{W}_1}$ is determined by a  minimal space-like surface $s\colon L^2\to \Sf_1^{n+2}$. 

The proof of the direct statement will be completed once we prove that, for each connected component $\mathcal{W}$ of the subsets $\mathcal{V}_1\cap \mathcal{W}^0_2\cap \mathcal{Y}^0_2$, $\mathcal{V}^0_2$, $\mathcal{V}_1\cap \mathcal{W}^0_2\cap \mathcal{U}^0_2$ and $\mathcal{V}_1\cap \mathcal{W}^0_2\cap \mathcal{Y}^0_1$, respectively, $f|_{\mathcal{W}}$ is a conformally surface-like hypersurface determined by a surface $h\colon L^2\to \Q_\epsilon$, $\epsilon \in \{-1,0,1\}$ of one of the following types:
   \begin{itemize}
   \item[\hypertarget{i}{(i)}]  an isothermic surface;
   \item[\hypertarget{ii}{(ii)}] a generalized cone over a unit-speed curve $\gamma\colon J\to \Q_{c}^2$ in an umbilical surface $\Q_{c}^2\subset \Q_{\epsilon}^3$, $c\geq \epsilon$;
   \item[\hypertarget{iii}{(iii)}]  a minimal surface;
   \item[\hypertarget{iv}{(iv)}]  an umbilical surface.
   \end{itemize}
Notice that any surface as in $(ii)$, $(iii)$ and $(iv)$ is also isothermic
(for $h$ as in $(ii)$ see Corollary $12$ in \cite{JT}). Notice also that $f|_{\mathcal{W}}$ being a conformally surface-like hypersurface determined by a surface $h\colon L^2\to \Q_\epsilon$, $\epsilon \in \{-1,0,1\}$, is equivalent to ${\mathcal{W}}$ being (isometric to)  a Riemannian product $L^2\times N^{n-2}$ and to $f|_{\mathcal{W}}$ being given by
$f|_{\mathcal{W}}=\mathcal{I}\circ\Phi\circ (h\times i)$, where $i$ is the inclusion map of an open subset $N^{n-2}$ of either $\R^{n-2}$ or $\mathbb{Q}_{-\epsilon}^{n-2}$, according to whether $\epsilon$ is zero or not, $\mathcal{I}$ is a Moebius transformation of $\R^{n+1}$,  and $\Phi$ is the standard isometry $\Phi\colon \R^3\times \R^{n-2}\to \R^{n+1}$ if $\epsilon=0$ and, if $\epsilon =-1$ or $1$, respectively,  the conformal diffeomorphism 
\begin{itemize}
\item $\Phi\colon \mathbb{H}^3\times \mathbb{S}^{n-2}\subset \mathbb{L}^{4}\times \mathbb{R}^{n-1}\to \mathbb{R}^{n+1}\setminus \R^{2}$, $\Phi(x, y)=\frac{1}{x_0}(x_1, x_2,y)
$
for all $x=x_0e_0+x_1e_1+x_2e_2+x_3e_3\in \mathbb{L}^{4}$ and $y=(y_1, \ldots, y_{n-1})\in \mathbb{S}^{n-2}\subset \mathbb{R}^{n-1}$, where $\{e_0, e_1, e_2, e_3\}$ is a pseudo-orthonormal basis of the Lorentzian space $\mathbb{L}^{k+1}$ with $\<e_0, e_0\>=0=\<e_3, e_3\>$ and $\<e_0, e_3\>=-1/2$.
\item $\Phi\colon \mathbb{S}^{3}\times \mathbb{H}^{n-2}\subset \mathbb{R}^{4}\times \mathbb{L}^{n-1}\to \mathbb{R}^{n+1}\setminus \R^{n-3}$, $\Phi(x, y)=\frac{1}{y_0}(x, y_0, \ldots, y_{n-3})$
for all $x=(x_1, \ldots, x_4)\in \mathbb{S}^3\subset \mathbb{R}^4$ and $y=y_0e_0+\ldots y_{n-3}e_{n-2}\in \mathbb{H}^{n-2}\subset\mathbb{L}^{n-1}$, where $\{e_0, \ldots,  e_{n-2}\}$ is a pseudo-orthonormal basis of $\mathbb{L}^{n-1}$ with $\<e_0, e_0\>=0=\<e_{n-2}, e_{n-2}\>$ and $\<e_0, e_{n-2}\>=~-1/2$.

\end{itemize}
\vspace{2ex}
   
\noindent \emph{Case $\hyperlink{i}{(i)}$:}  Let $\mathcal{W}$ be a connected component of
$\mathcal{V}_1\cap \mathcal{W}^0_2\cap \mathcal{Y}^0_2$.  Since, in particular, 
$\mathcal{W}\subset \mathcal{W}_2$, then $C(\Gamma(\Delta))\subset \spa\{I\}$ on $\mathcal{W}$.  By Lemma~\ref{le:split}, $f|_{\mathcal{W}}$ is a conformally surface-like hypersurface determined by a surface $h\colon L^2\to \Q_\epsilon$, $\epsilon \in \{-1,0,1\}$.
Thus, with notations as in the preceding paragraph,  we may write $\mathcal{W}= L^2\times N^{n-2}$ and $f|_{\mathcal{W}}=\mathcal{I}\circ\Phi\circ (h\times i)$. In particular, the distributions $\Delta$ and  ${\Delta}^\perp$ are given by the tangent spaces to  $N^{n-2}$ and $L^2$, respectively. 

Denote by $g_1$ the product metric of $\mathbb{Q}_\epsilon^3\times \mathbb{Q}_{-\epsilon}^{n-2}$ and let $g_2$ be the metric on $\mathbb{Q}_\epsilon^3\times \mathbb{Q}_{-\epsilon}^{n-2}$ induced from the metric of $\R^{n+1}$ by the conformal diffeomorphism $\mathcal{I}\circ \Phi$. Let $\varphi \in C^{\infty}(\mathbb{Q}_\epsilon^3\times \mathbb{Q}_{-\epsilon}^{n-2})$
be the conformal factor of $g_2$ with respect to $g_1$, that is, $g_2=\varphi^2 g_1$. Then the shape operators of  $F_1=h\times i\colon L^2\times N^{n-2}\to \mathbb{Q}_\epsilon^3\times \mathbb{Q}_{-\epsilon}^{n-2}$ and of
$F_2=F_1\colon (L^2\times N^{n-2}, F_1^*g_2)\to (\mathbb{Q}_\epsilon^3\times \mathbb{Q}_{-\epsilon}^{n-2}, g_2)$ with respect to unit normal vector fields $N_1$ and $N_2=N_1/\varphi$, respectively, are
related by
\be\label{eq:relsff}
A_{N_2}^{F_2}=\frac{1}{\varphi\circ F_1}A^{F_1}_{N_1}-\frac{g_1(\nabla^1 \varphi, N_1)}{(\varphi\circ F_1)^2} I.
\ee
We recall that the Levi-Civita connections $\bar\nabla$ and $\nabla$ of the metrics $\bar g$ and $g=\<\,,\,\>$ on $M^n=L^2\times N^{n-2}$ induced by $F_1$ and $F_2$, respectively, satisfy
\begin{equation}\label{eq:conn}
\nabla_XY=\bar\nabla_XY+\frac{1}{\bar\varphi}(X(\bar\varphi)Y+Y(\bar\varphi)X-\bar g(X,Y)\bar\nabla\bar\varphi),
\end{equation}
where $\bar\varphi=\varphi\circ F_1$. It follows from \eqref{eq:conn} that the mean curvature vector field $\delta\in\Gamma(\Delta^\perp)$ of $\Delta$ (with respect to the metric $g$) is
\be\label{eq:delta}
\delta=-\bar{\varphi}^{-3}(\bar\nabla\bar\varphi)_{\Delta^\perp}.
\ee

     Now, we can write \eqref{hyp2} as
\begin{equation} \label{eq:nb}
(\nabla_X \bar{\mathcal{B}})Y-(\nabla_Y\bar{\mathcal{B}})X+(X\wedge Y)(A\nabla\rho-\nabla b)=0,
\end{equation}
for all $X,Y\in\mathfrak{X}(M)$. The $\Delta$-component of \eqref{eq:nb} evaluated in unit vector fields $Z\in\Gamma(\Delta^\perp)$ and $T\in\Gamma(\Delta)$ gives
$
\<\bar{\mathcal{B}}Z,\nabla_TT\>=\<Z,A\nabla\rho-\nabla b\>,
$
or equivalently,
\be\label{eq:bdelta}
\bar{\mathcal{B}}\delta=A\nabla\rho-\nabla b.
\ee

 Since $\mathcal{W}\subset \mathcal{Y}^0_2$, there exists locally a smooth function $\theta$ and an orthonormal frame $\{X, Y\}$ of $\Delta^\perp$ given by principal directions of  
$f$ such that $\bar{\mathcal{B}} X=\theta Y$ and $\bar{\mathcal{B}} Y=\theta X$.
From \eqref{eq:conn} we have $\<\nabla_XT,X\>=\<\nabla_YT,Y\>=T(\log\circ\bar{\varphi})$.
Evaluating \eqref{eq:nb} in $T$ and $X$ (or $Y$) gives
\be\label{eq:ttheta}
T(\theta)=-T(\log\circ\bar{\varphi})\theta,
\ee
whereas \eqref{eq:nb} evaluated in $X$ and $Y$ yields
\be\label{eq:xtheta}
X(\theta)=2\theta\<\nabla_YY,X\>-\<Y,A\nabla\rho-\nabla b\>
\ee
and
\be\label{eq:ytheta}
Y(\theta)=2\theta\<\nabla_XX,Y\>-\<X,A\nabla\rho-\nabla b\>.
\ee
Set $\bar{\theta}=\bar{\varphi}\theta$ . It follows from \eqref{eq:ttheta} that 
$T(\bar{\theta})=0$ for any $T\in\Gamma(\Delta)$.
Thus $\bar{\theta}$ induces a function on $L^2$, which we also denote by $\bar \theta$.
By \eqref{eq:relsff}, the vector fields $\bar{X}=\bar{\varphi}X$ and $\bar{Y}=\bar{\varphi}Y$ form an orthonormal frame of principal directions of $h$. Using \eqref{eq:conn}, \eqref{eq:delta}, \eqref{eq:bdelta} and  \eqref{eq:xtheta} we obtain
\begin{align}\label{eq:barthetcod1}
    \bar{X}(\bar{\theta})&=\bar{X}(\bar{\varphi})\theta+\bar{\varphi}^2X(\theta)\nonumber\\
    &=\bar{X}(\bar{\varphi})\theta+\bar{\varphi}^2(2\theta\<\nabla_YY,X\>-\<Y,A\nabla\rho-\nabla b\>)\nonumber\\
    &=\bar{X}(\bar{\varphi})\theta+\bar{\varphi}^2(2\theta(\bar{\varphi}^{-1}\bar{g}(\bar{\nabla}_{\bar{Y}}\bar{Y},\bar{X})-\bar{\varphi}^{-2}\bar{X}(\bar{\varphi}))-\<\bar{B}Y,\delta\>)\nonumber\\
    &=\bar{X}(\bar{\varphi})\theta+2\bar{\theta}\bar{g}(\bar{\nabla}_{\bar{Y}}\bar{Y},\bar{X})-2\theta\bar{X}(\bar{\varphi})-\bar{\varphi}^2\theta\<X,\delta\>\nonumber\\
    &=\bar{X}(\bar{\varphi})\theta+2\bar{\theta}\bar{g}(\bar{\nabla}_{\bar{Y}}\bar{Y},\bar{X})-2\theta\bar{X}(\bar{\varphi})+\theta\bar{\varphi}^{-1}\<X,\bar{\nabla}\bar{\varphi}\>\nonumber\\
    &=\bar{X}(\bar{\varphi})\theta+2\bar{\theta}\bar{g}(\bar{\nabla}_{\bar{Y}}\bar{Y},\bar{X})-2\theta\bar{X}(\bar{\varphi})+\theta\bar{X}(\bar{\varphi})\nonumber\\
    &=2\bar{\theta}\bar{g}(\bar{\nabla}_{\bar{Y}}\bar{Y},\bar{X}).
\end{align}
 A similar computation  using \eqref{eq:ytheta} instead of  \eqref{eq:xtheta} gives
\be\label{eq:barthetcod2}
\bar{Y}(\bar{\theta})=2\bar{\theta}\bar{g}(\bar{\nabla}_{\bar{X}}\bar{X},\bar{Y}).
\ee

Let $\mathcal{B}^*\in\Gamma(\End(TL))$ be defined by $\mathcal{B}^* X=\bar{\theta}\bar{Y}$ and $\mathcal{B}^* Y=\bar{\theta}\bar{X}$.
Then \eqref{eq:barthetcod1} and \eqref{eq:barthetcod2} are equivalent to $\mathcal{B}^*$ being a Codazzi tensor on $L^2$.
By $\eqref{eq:relsff}$, the shape operator $A^h$ of $h$ is a multiple of $\bar{A}|_{\Delta^\perp}$,  which has rank two, for $\mathcal{W}\subset \mathcal{V}_1$, and  $\mathcal{B}^*$ is a multiple of $\bar{\mathcal{B}}|_{\Delta^\perp}$. Therefore, the fact that $\bar{A}$ and $\bar{\mathcal{B}}$ satisfy \eqref{hyp1b} implies that  $\mathcal{B}^*$ and $A^h$  also satisfy \eqref{hyp1b} with respect to $\bar{g}$.  Since, in addition, $\trace \mathcal{B}^*=0$, it follows from  Proposition $8$ of \cite{JT}, together with Theorem $4.7$ in \cite{DJ1} (also stated  in \cite{JT} as Theorem $2$), that $h$ is locally infinitesimally Bonnet bendable,  hence isothermic by Proposition~$9$ of \cite{JT}.\vspace{2ex}

\noindent \emph{Case $\hyperlink{ii}{(ii)}$:}  First we show that the interior $\mathcal{V}^0_2$ of the subset $\mathcal{V}_2$  where $\dim \ker \bar{A}=n-1$ is contained in $\mathcal{W}^0_2\cap  \mathcal{Y}^0_2$. 
Clearly,  $\mathcal{V}_2\subset \mathcal{U}_1$, the subset where $\trace \bar{B}=0$, for if $\trace \bar{A}(x)=0$, that is, if $x\in \mathcal{V}_2\cap \mathcal{U}_2$, then we would have $\bar{A}(x)=0$, in contradiction with the assumption that $f$ is free of umbilic points. Also, since $\dim \ker \bar{A}=n-1$ on $\mathcal{V}_2$, then $\bar{A}|_{\Delta^\perp}$ can not be a multiple of the identity endomorphism $I\in \Gamma(\End(\Delta^\perp))$ at any point of $\mathcal{V}_2$. Thus $\mathcal{V}_2\subset \mathcal{Y}_2$. 

  We now show that, on any connected component  $\mathcal{W}$ of $\mathcal{V}^0_2$, the splitting tensor $C\colon\Gamma(\Delta)\to\Gamma(\End(\Delta^\perp))$ of $\Delta$ satisfies 
  $C(\Gamma(\Delta))\subset \spa\{I\}$, which will yield the inclusion $\mathcal{V}^0_2\subset \mathcal{W}^0_2$.  So let $\mathcal{W}$ be such a connected component. Since we already know that $\mathcal{W}\subset \mathcal{Y}^0_2$, there exist locally smooth functions $\theta, \mu$ and  unit vector fields $Y\in \Gamma(\Delta^\perp\cap \ker \bar{A})$ and $X\in \Gamma((\ker{\bar A})^\perp)$ such that $\bar{A}X=\mu X$, 
$\bar{\mathcal{B}}X=\theta Y$ and $\bar{\mathcal{B}}Y=\theta X$.

   Applying \eqref{eq:nb} to $T$ and $Y$ gives
\begin{equation} \label{eq:yt}
\begin{array}{l}T(\theta)X+\theta\nabla_TX-\bar{\mathcal{B}}\nabla_TY+\bar{\mathcal{B}}\nabla_YT+\<\lambda\nabla\rho-\nabla b, Y\>T\vspace{1ex}\\
\hspace*{20ex}-\<\lambda\nabla\rho-\nabla b, T\>Y=0.
\end{array}
\end{equation}
Taking the $X$-component of \eqref{eq:yt} we obtain
\begin{equation}\label{eq:nyyt}
T(\theta)=\theta\<\nabla_YY, T\>,
\end{equation}
whereas the $Y$-component and the $T$-component give, respectively,
\begin{equation}\label{eq:lnrt}
\<\lambda\nabla\rho-\nabla b, T\>=0
\end{equation}
and
$$
\<\lambda\nabla\rho-\nabla b, Y\>=\theta\<\nabla_TT, X\>.
$$

Applying \eqref{eq:nb} to $T$ and $X$ and using \eqref{eq:lnrt} give
\begin{equation} \label{eq:xt}
\begin{array}{l}T(\theta)Y+\theta\nabla_TY-\bar{\mathcal{B}}\nabla_TX
+\bar{\mathcal{B}}\nabla_XT+\<\mu\nabla\rho-\nabla b, X\>T=0.
\end{array}
\end{equation}
Taking the $S$-component of \eqref{eq:xt} for $S\in \Gamma(\Delta)$ with $\<S, T\>=0$ gives
$$\<\nabla_T S, Y\>=0,$$
and taking its $T$-component yields
$$
\theta\<\nabla_TT, Y\>=\<\mu\nabla\rho-\nabla b, X\>.
$$
Using that $\ker\bar{A}=\{X\}^\perp$ is an umbilical distribution, it follows that the same holds for $\Delta$.

Taking the $X$-component of \eqref{eq:xt} yields
\begin{equation}\label{eq:nxyt}
\<\nabla_XY, T\>=0,
\end{equation}
whereas the $Y$-component gives
\begin{equation}\label{eq:nxxt}
T(\theta)=\theta\<\nabla_XX, T\>.
\end{equation}
It follows from \eqref{eq:nyyt},  \eqref{eq:nxyt} and  \eqref{eq:nxxt}, taking into account that one also has $\<\nabla_YX, T\>~=~0$, that
the distribution $\Delta^\perp$ is umbilical with mean curvature vector field 
$\zeta=(\nabla \log \theta)|_{\Delta}$, which is equivalent to the splitting tensor $C\colon\Gamma(\Delta)\to\Gamma(\End(\Delta^\perp))$ of $\Delta$ satisfying $C_T=\<\zeta, T\> I$ for all $T\in\Gamma(\Delta)$. 

Now that we know that $\mathcal{V}^0_2\subset \mathcal{W}^0_2\cap  \mathcal{Y}^0_2$, the argument used in case $\hyperlink{i}{(i)}$ shows that $f|_{\mathcal{W}}$ is a conformally surface-like hypersurface determined by an isothermic surface $h\colon L^2\to \Q^3_\epsilon$, $\epsilon \in \{-1,0,1\}$.
But since $\rank \ker \bar{A}=n-1$ on $\mathcal{V}_2$, then $h$ has index of relative nullity equal to one at any point. By Corollary~$12$ in \cite{JT}, $h$ is a generalized cone over a unit-speed curve $\gamma\colon J\to \Q_{c}^2$ in an umbilical surface $\Q_{c}^2\subset \Q_{\epsilon}^3$, $c\geq \epsilon$. \vspace{2ex}

\noindent \emph{Cases $\hyperlink{iii}{(iii)}$ and $\hyperlink{iv}{(iv)}$:} Let $\mathcal{W}$ be a connected component of $\mathcal{V}_1\cap \mathcal{W}^0_2\cap \mathcal{U}_2^0$ (respectively, $\mathcal{V}_1\cap \mathcal{W}^0_2\cap \mathcal{Y}^0_1$). As in cases $\hyperlink{i}{(i)}$ and $\hyperlink{ii}{(ii)}$, $f|_{\mathcal{W}}$ is a conformally surface-like hypersurface determined by a surface $h\colon L^2\to \Q_\epsilon$, $\epsilon \in \{-1,0,1\}$, by Lemma \ref{le:split}.  Since $\trace \bar{A}=0$ on $\mathcal{U}_2$ (respectively, $\bar{A}|_{\Delta^\perp}$ is a multiple of the identity endomorphism of $\Delta^\perp$ on $\mathcal{Y}_1$), it follows from \eqref{eq:relsff} that also $\trace A^h=0$ (respectively, $A^h$ is  a multiple of the identity endomorphism of $TL$), hence $h$ is a minimal surface (respectively, $h$ is umbilical). \vspace{1ex}
 
  We now prove the converse. Assume first that $f\colon M^n\to\R^{n+1}$ is a simply connected hypersurface whose central sphere congruence is determined by a minimal space-like surface $s\colon L^2\to \Sf_{1,1}^{n+2}$.
Let $\bar{J}\in\Gamma(\End(TL))$ represent a rotation of angle $\pi/2$, and let $\bar{X},\bar{Y}$ be an orthonormal frame satisfying $\bar{J}\bar{X}=\bar{Y}$ and $\bar{J}\bar{Y}=-\bar{X}$. Then $\bar{J}$ is parallel with respect to the  Levi-Civita connection  $\nabla'$ on $L^2$,  
hence it is, in particular, a Codazzi tensor on $L^2$.
Since $s$ is minimal, then
$
\a'(\bar{X},\bar{X})+\a'(\bar{Y},\bar{Y})=0,
$
hence $s$ is  a special elliptic surface by Proposition $11$ in \cite{DJV}.

By Theorem $1$ in \cite{DJV},  $f$ admits a non-trivial conformal infinitesimal bending $\T$. We now show that $\T$ is also an infinitesimal Moebius bending.
Let $X,Y\in\mathfrak{X}(M)$ be the lifts of $\bar{X}$ and $\bar{Y}$. From \eqref{relmetr} we see that $\bar{A}X$ and $\bar{A}Y$ form an orthonormal frame of $\Delta^\perp$.
Let $J\in\Gamma(\End(\Delta^\perp))$ be the lift of $\bar{J}$. It was shown in the proof of the converse of Theorem $1$ in \cite{DJV} that $\bar{A}J$ is symmetric. Thus
$$\<J\bar{A}X,\bar{A}X\>=\<\bar{A}J\bar{A}X,X\>=\<\bar{A}X,\bar{A}JX\>=\<\bar{X},\bar{J}\bar{X}\>'=0.$$
Similarly,  $\<J\bar{A}Y,\bar{A}Y\>=0$, $\<J\bar{A}X,\bar{A}Y\>=-1$ and $\<J\bar{A}Y,\bar{A}X\>=1$. Hence $J$ is an orthogonal tensor, and the symmetry of $\bar{A}J$ implies that $\trace\bar{A}=0$. By Proposition \ref{prop:imbhyp},  $\T$ is an infinitesimal Moebius bending.
 
    Now let $f\colon M^n\to \R^{n+1}$ be a conformally surface-like hypersurface determined by an isothermic surface $h\colon L^2\to \Q_\epsilon$, $\epsilon \in \{-1,0,1\}$.
Then $h$ is locally infinitesimally Bonnet bendable (see, e.g., Proposition $9$ and Remark $11$ of \cite{JT}), that is, it admits locally a non-trivial infinitesimal variation $h_t\colon L^2\to \Q_\epsilon$ such that the metrics $\bar{g}_t$ induced by the ${h_t}'s$ and their mean curvatures $\mathcal{H}_{t}$  satisfy 
$
\partial/\partial t|_{t=0}\bar{g}_t=0=\partial/\partial t|_{t=0}\mathcal{H}_t.
$
Thus also $\partial/\partial t|_{t=0}K_t=0$, where $K_t$ is the Gauss curvature of $\bar{g}_t$. 

     Let $f_t$ be the variation of $f$ given by the conformally surface-like hypersurfaces determined by $h_t$. The Moebius metric of $f_t$ is  (see Remark $3.7$ in \cite{LMW})
$$
\left(4\mathcal{H}_t^2-\frac{2n}{n-1}(K_t-\epsilon)\right)(\bar{g}_t+g_{-\epsilon}).
$$
Here $g_{-\epsilon}$ is the metric of $\Q_{-\epsilon}^{n-2}$ and $\bar{g}_t+g_{-\epsilon}$ denotes the product metric on $L^2\times \Q_{-\epsilon}^{n-2}$.
Therefore the immersions $f_t$ determine an infinitesimal Moebius variation of $f$. It remains to argue that the latter is non-trivial.

From Proposition \ref{eq:alinha} we know that the associated tensor $\mathcal{B}$  satisfies \eqref{BAlinha}. On the other hand, by \eqref{eq:relsff} the shape operator $A_t$ of $f_t$ has the form 
\begin{equation}\label{eq:at}
A_t=\delta_{1}(t)\bar{A}_t+\delta_2(t)I,
\end{equation}
for some smooth functions $\delta_1$ and $\delta_2$, with $\delta_1(0)\neq 0$. Here $\bar{A}_t$ denotes the second fundamental form of $h_t$  extended to $TM$ by defining $\bar{A}_tT=0$ for any $T$ tangent to $\Q_{-\epsilon}^{n-2}$.
Since $\partial/\partial t|_{t=0} \bar{A}_t\neq 0$, for $h_t$ determine a non-trivial infinitesimal Bonnet variation of $h$, it follows from \eqref{BAlinha} and \eqref{eq:at} that 
$\mathcal{B}$ is not a multiple of the identity endomorphism. Hence the infinitesimal Moebius variation of $f$ determined by $f_t$ is non-trivial (see Remarks \ref{re:infvar}--$2)$).
\qed

\vspace*{5ex}

\noindent Universidade de S\~ao Paulo\\
Instituto de Ci\^encias Matem\'aticas e de Computa\c c\~ao.\\
Av. Trabalhador S\~ao Carlense 400\\
13566-590 -- S\~ao Carlos\\
BRAZIL\\
\texttt{mibieta@icmc.usp.br} and \texttt{tojeiro@icmc.usp.br}

\newpage


\begin{thebibliography}{lll}

\bibitem{ca1} Cartan, E., 
\emph{La d\'eformation des hypersurfaces dans l'espace euclidien 
r\'eel a $n$ dimensions},
Bull. Soc. Math. France 44 (1916), 65--99.

\bibitem{Ca2} Cartan, E., 
\emph{La d\'eformation des hypersurfaces dans l'espace conforme r\'eel  
a $n\ge 5$ dimensions}, 
Bull. Soc. Math. France 45 (1917), 57--121.

\bibitem{DJ} Dajczer M. and Jimenez M. I.,
\emph{Conformal infinitesimal variations of 
submanifolds}, Differ. Geom. Appl. (2021),

\bibitem{DJ1} Dajczer M. and Jimenez M. I.,
\emph{Infinitesimal variations of 
submanifolds}, Ensaios Matemáticos, 35 (2021), 1-156.

\bibitem{DJV} Dajczer M., Jimenez M. I. and Vlachos, Th.,
\emph{Conformal infinitesimal variations of Euclidean hypersurfaces},
 Ann. Mat. Pura Appl. 201 (2022), 743 -768.

\bibitem{DV} Dajczer M. and Vlachos, Th.,
\emph{Infinitesimally bendable Euclidean hypersurfaces}, Ann. Mat. Pura Appl. 196 (2017), 1961--1979 and Ann. Mat. Pura Appl. 196 (2017), 1981--1982.

\bibitem{DFT}  Dajczer, M., Florit, L. and  Tojeiro, R., 
 \emph{On a class of submanifolds carrying an
extrinsic umbilic foliation\/}. Israel J. Math. 125 (2001),
203-220.

\bibitem{DFT2} Dajczer, M., Florit, L. and Tojeiro, R.,
\emph{On deformable hypersurfaces in space forms},
Ann. Mat. Pura Appl.  174 (1998), 361--390.

\bibitem{DT} Dajczer, M. and  Tojeiro, R., 
Submanifold theory beyond an introduction.
Universitext. Springer, New York, 2019.

\bibitem{DT1} Dajczer, M. and Tojeiro, R.,
\emph{On Cartan's conformally deformable hypersurfaces},
Michigan Math. J.  47 (2000), 529--557.


\bibitem{JT}  Jimenez, M. I. and  Tojeiro, R., 
 \emph{Infinitesimally Bonnet bendable hypersurfaces\/}. J. Geom. Anal. 33, no. 5, article 140 (2023). 
 
 \bibitem{JT2}  Jimenez, M. I. and  Tojeiro, R., 
 \emph{On the Moebius deformable hypersurfaces\/}. To appear in Rev. Mat. Iberoamericana.

\bibitem{LMW} Li, T., Ma, X. and Wang, C., 
\emph{Deformations of hypersurfaces preserving the M\"obius metric and
a reduction theorem},  
Adv.  Math. 256 (2014), 156--205.

\bibitem{sb0} Sbrana, U.,
\emph{Sulla deformazione infinitesima delle ipersuperficie},
Ann. Mat. Pura  Appl. 15 (1908), 329--348.

\bibitem{sb} Sbrana, U.,
\emph{Sulle variet\`a ad $n-1$ dimensioni deformabili nello 
spazio euclideo ad $n$ dimensioni},
Rend. Circ. Mat. Palermo  27 (1909), 1--45.

\bibitem{Wa}  Wang, C.,  
\emph{M\"obius geometry of submanifolds in $S^n$},  
Manuscripta Math. 96 (1998), 517--534.
\end{thebibliography}
\end{document}